# ISOMORPHISMS OF POISSON SYSTEMS OVER LOCALLY COMPACT GROUPS

BY AMANDA WILKENS[a]

*For my brother, Daniel, and in memory of him.*

*Mathematics, University of Texas at Austin*, [a]*amanda.wilkens@math.utexas.edu*

A Poisson system is a Poisson point process and a group action, together forming a measure-preserving dynamical system. Ornstein and Weiss proved Poisson systems over many amenable groups were isomorphic in their 1987 paper. We consider Poisson systems over nondiscrete, noncompact, locally compact Polish groups, and we prove by construction all Poisson systems over such a group are finitarily isomorphic, producing examples of isomorphisms for nonamenable group actions. As a corollary, we prove Poisson systems and products of Poisson systems are finitarily isomorphic.

For a Poisson system over a group belonging to a slightly more restrictive class than above, we further prove it splits into two Poisson systems whose intensities sum to the intensity of the original, generalizing the same result for Poisson systems over Euclidean space proved by Holroyd, Lyons and Soo.

**1. Introduction.** In what follows, we construct an isomorphism between any two Poisson systems over a nondiscrete, noncompact, locally compact Polish group. We define Poisson systems and isomorphisms between them before formally stating our main result. Let $\gamma > 0$. A random variable $Z$ taking values on $\mathbb{N} := \{0, 1, 2, \ldots\}$ with

$$\mathbb{P}(Z = m) = \frac{e^{-\gamma}\gamma^m}{m!}$$

is a *Poisson random variable with mean $\gamma$*. Let $G$ be a group with $\sigma$-finite Haar measure $\lambda$, and let $\alpha > 0$. A *Poisson point process on $G$ with intensity $\alpha$* is a random variable $X$, which takes values on the space $\mathbb{M} := \mathbb{M}(G, \mathcal{B}(G))$ of Borel simple point measures on $G$, such that: For every Borel subset $A \in \mathcal{B}(G)$ with $\lambda(A) < \infty$, the number of points of $X$ in $A$, denoted by $N(A)$, is a Poisson random variable with mean $\alpha \cdot \lambda(A)$. For pairwise disjoint Borel sets $A_1, \ldots, A_\ell \in \mathcal{B}(G)$, the random variables $N(A_1), \ldots, N(A_\ell)$ are independent.

Let $P_\alpha$ be the law of a Poisson point process on $G$ with intensity $\alpha$, and let $\mu \in \mathbb{M}$ and $g \in G$. The triple $(\mathbb{M}, P_\alpha, G)$ is a measure-preserving dynamical system, which we refer to as a *Poisson system* and denote as $\mathbb{X}_\alpha$, where the action of $G$ on $\mathbb{M}$ is given by $g(\mu)(A) = \mu(g^{-1}A)$ for all $A \in \mathcal{B}(G)$. It may be helpful to think of a Poisson system as a random collection of points in a state space, along with actions inherited from the state space, which left-translate the collection of points around the space.

A map $\phi : \mathbb{X}_\alpha \to \mathbb{X}_\beta$ is a *factor from $\mathbb{X}_\alpha$ to $\mathbb{X}_\beta$* if, on a set of $P_\alpha$ full-measure, it behaves with respect to the group action and relevant measures, so that $\phi \circ g = g \circ \phi$ and $P_\alpha \circ \phi^{-1} = P_\beta$. A factor is an *isomorphism* if it is a bijection and its inverse is a factor from $\mathbb{X}_\beta$ to $\mathbb{X}_\alpha$.

We impose some additional restrictions on $G$ and, by the end of Section 2, will have detailed each choice. Let $G$ be nondiscrete, noncompact and locally compact Polish.

---

THEOREM 1.1. *All Poisson systems over G are isomorphic.*

Ornstein and Weiss proved all Poisson systems over any group with "good entropy theory" were isomorphic [17]. Groups with good entropy theory are locally compact, amenable and unimodular, with some further restrictions. There are many new isomorphism examples as a result of Theorem 1.1; some in particular are for Poisson systems over the Euclidean group $E(n)$ for $n \geq 3$ and the special linear group $\mathrm{SL}(n, \mathbb{R})$ for $n \geq 2$. Any real Lie group with countably many connected components satisfies the conditions on $G$; so does the field $\mathbb{Q}_p$ of $p$-adic numbers for any prime $p$. Our construction does not distinguish between amenable and nonamenable groups, or unimodular and nonunimodular.

We prove the next theorem over the course of proving Theorem 1.1.

THEOREM 1.2. *Let $\alpha, \beta > 0$. There exists a factor from $\mathbb{X}_\alpha$ to $\mathbb{X}_\beta$.*

Holroyd, Lyons and Soo construct a factor from $\mathbb{X}_\alpha$ to $\mathbb{X}_\beta$ in the case $G = \mathbb{R}^n$, for all $n \geq 1$, in their paper on Poisson splitting [8]. The Ornstein and Weiss isomorphisms of [17], while applicable to a large class of groups, are nonconstructive. Via an explicit construction, Holroyd, Lyons and Soo gain some additional properties. The factor maps of [8] extend the action from translations to isometries of $\mathbb{R}^n$ and are finitary, which we define shortly. Our proof of Theorem 1.2 extends the construction in [8]. In [22], Soo and the author construct an isomorphism between any two Poisson systems over $\mathbb{R}^n$. The isomorphism builds on the factor map from [8] and carries its additional properties. Our proof of Theorem 1.1 extends this construction.

Recall a locally compact space is Polish if and only if it is second countable. Birkhoff and Kakutani proved a topological group is metrizable if and only if it is first countable ([5], Theorem 2.B.2). A metric $d$ on a topological space $(X, T)$ is *proper* if balls with respect to $d$ are relatively compact, and is *compatible* if $d$ defines $T$. Struble proved a locally compact topological group is metrizable with a proper, left-invariant, compatible metric if and only if it is second countable [24].

To define a finitary factor, we fix $d$ as a proper, left-invariant, compatible metric on $G$ and denote the open ball centered at $g \in G$ of radius $r$, with respect to the metric $d$, as $B(g, r)$ (although, for most of the paper, we use a more specific metric and a different definition of $B(g, r)$, beginning in Section 3). We use either the random variable $X$, which takes values on $\mathbb{M}$, or the collection of Borel simple point measures $\mu \in \mathbb{M}$ to refer to a Poisson point process on $G$, depending on context. The restriction of $\mu \in \mathbb{M}$ to a Borel set $A$ is $\mu|_A(B) := \mu(B \cap A)$ for $B \in \mathcal{B}(G)$ and is itself a Poisson process on $A$ whenever $\mu$ is a Poisson process on $G$ ([15], Chapter 2).

Let $\phi$ be a factor from $\mathbb{X}_\alpha$ to $\mathbb{X}_\beta$, and let $\mu, \mu' \in \mathbb{M}$. We always denote the identity element of a group as $e$. A *coding window* of $\phi$ is a map $\omega : \mathbb{M} \to \mathbb{N} \cup \{\infty\}$ such that if

$$\mu|_{B(e, \omega(\mu))} = \mu'|_{B(e, \omega(\mu))}$$

then $\phi(\mu)|_{B(e,1)} = \phi(\mu')|_{B(e,1)}$. The factor map $\phi$ is *d-finitary* if there exists a $P_\alpha$ almost surely finite coding window of $\phi$, or in other words, if $P_\alpha$ almost surely the ball $B(e, 1)$ under $\mathbb{X}_\beta$ is completely determined by a ball at $e$ of finite radius under $\mathbb{X}_\alpha$, with respect to the metric $d$. The choice of $e$ is arbitrary due to $G$-equivariance of $\phi$. Any two proper, left-invariant, compatible metrics on $G$ are coarsely equivalent; in particular, if $\phi$ is a d-finitary factor with respect to one such metric, it is also a d-finitary factor with respect to the other ([5], Definition 3.A.3 and Corollary 4.A.6). We say a factor is *finitary* if it is $d$-finitary with respect to any left-invariant, proper, compatible metric $d$ (so $d$-finitary immediately implies finitary). A finitary factor is a finitary isomorphism if it is a bijection and its inverse is a finitary factor from $\mathbb{X}_\beta$ to $\mathbb{X}_\alpha$.

Theorem 1.3 is a direct consequence of the proofs of Theorems 1.1 and 1.2, wherein we extend the construction of the finitary maps between Poisson systems over $\mathbb{R}^n$ from [8, 22].

THEOREM 1.3. *All Poisson systems over G are finitarily isomorphic.*

1.1. *The isomorphism problem for countable groups.* Theorems 1.1, 1.2 and 1.3 have a rich history in the setting of countable groups, where Poisson systems take the form of Bernoulli shifts. Let $D$ be a countable group and $(K, \kappa)$ a standard Borel probability space. Endow $K^D := \{x : D \to K\}$ with the product $\sigma$-algebra generated by cylinder sets, and let $\kappa^D$ be the product measure on $K^D$. The group $D$ acts on $K^D$ by $(dx)(g) = x(d^{-1}g)$ for all $x \in K^D$ and $d, g \in D$. The system $(K^D, \kappa^D, D)$ is called a *Bernoulli shift over D*. It was an open question as to whether the Bernoulli shifts over the integers with uniform outcome sets $\{0, 1\}$ and $\{0, 1, 2\}$ were isomorphic as measure-preserving dynamical systems for about twenty years, until the late 1950s.

Around that time, Kolmogorov formalized the concept of entropy for ergodic theory [9]. Entropy measures the randomness of a system and is invariant under isomorphism in Kolmogorov's construction. Entropy does not apply to the isomorphism question in our setting, so we do not give a formal definition. But we note the entropies of the Bernoulli shifts referenced above are different, implying the shifts are not isomorphic. (Poisson systems over groups with good entropy theory have infinite entropy [17], and one might conjecture the same holds for Poisson systems over $G$.) Kolmogorov further proved entropy is nonincreasing under factor maps between Bernoulli shifts over $\mathbb{Z}$.

More results followed, including Sinai's on the existence of a factor map between shifts with nonincreasing entropy [9], and Ornstein's proof that shifts are isomorphic if and only if they share the same entropy [16]. Keane and Smorodinsky later constructed isomorphisms between shifts with finite coding windows [10, 11]. The finitary definition for Poisson systems first appears in [8] under the terminology "strongly finitary" and is based on the finitary definition for Bernoulli shifts, which can be found in [19].

In 1975, Kieffer proved entropy is nonincreasing under factor maps between shifts over any countable amenable group [14], and Stepin introduced a method called coinduction, which we use as well (see Section 5.2), to extend Ornstein's result to all countable groups containing $\mathbb{Z}$ as a subgroup and some further groups [23]. In their 1987 paper, Ornstein and Weiss, as a precursor to proving Poisson systems over any group with good entropy theory were isomorphic, proved Bernoulli shifts with equal entropy over such a group were isomorphic [17].

More recently, Bowen [2] and Kerr and Li [13] extended entropy theory to countable sofic groups. Bowen additionally proved shifts with equal entropy over any countable group, as long as both shifts have outcome sets of cardinality at least 3, were isomorphic [4]. Seward removed the cardinality restriction and added the finitary property in [21], thus closing the open question for the countable sofic setting. In [20], Seward states a conjecture which, if true, closes the open question in full generality for countable groups.

1.2. *Paper outline.* We first push the factor map and isomorphism constructions on $\mathbb{R}^n$ in [8, 22] through to nondiscrete, noncompact and locally compact Polish groups that are compactly generated (Section 5.1); see Sections 2 and 3 for more information on the compactly generated requirement. From here, we drop the compactly generated condition with a coinduction argument for the case in which $G$ contains a noncompact, compactly generated open subgroup (Section 5.2), and a modified approach for the case in which $G$ does not (Section 6.1). We discuss technical details of Poisson point processes and systems on groups, obstructions to lifting any requirements on $G$ and a further question in Section 2.

We begin construction of an isomorphism between Poisson systems over compactly generated $G$ in Section 3 by building a $G$-equivariant partition deterministically from information in the source system. The idea goes back to Keane and Smorodinsky's finitary factor maps for Bernoulli shifts in the $\mathbb{Z}$ setting, as found in [10]. Here, Keane and Smorodinsky introduce markers and fillers; markers are strings of symbols defined to identify potential sources of randomness, and fillers occupy the remaining space. To be useful, markers must be defined carefully enough as to not disrupt inherent independence conditions of the source system. We prove independence properties and harness the randomness made available by our marker system in Section 4.

To convert pieces of randomness into uniform random variables, we rely on the Borel isomorphism theorem for measures (see Remark 2.1) and the uncountable property of $G$ (since $G$ is nondiscrete and second countable, it must be uncountable). We apply the isomorphism theorem again to convert the uniform random variables into Poisson point processes on bounded Borel sets determined by the marker partition. As long as we are careful with independence properties, the processes sum to a Poisson point process on $G$ (see Remark 3.5). This is the idea behind the proof of Theorem 1.2.

With Theorems 1.2, 1.3 and only a bit more work, we extend the Poisson splitting results in [8] from $\mathbb{R}^n$ to nondiscrete, noncompact, locally compact Polish groups that are either compactly generated or contain a noncompact, compactly generated open subgroup. A *Poisson splitting* is a map with the properties described in Theorem 1.4 (but not necessarily finitary or a factor). For example, if one takes a Poisson point process on $G$ with intensity $\alpha$, and for some $\beta < \alpha$, colors a point red with probability $\beta/\alpha$ and blue with probability $1 - \beta/\alpha$ by means of randomness independent of the process, then the map sending the original process to just the red (or blue) points is a Poisson splitting. For a more detailed picture of these maps, see [8].

We discuss the extra group conditions, prove Theorem 1.4 and ask whether the extra conditions are necessary in Section 7.

THEOREM 1.4. *Let $G$ be a nondiscrete, noncompact, locally compact Polish group, either compactly generated or containing a noncompact, compactly generated open subgroup. For all $\alpha > \beta > 0$, there exists a monotone, finitary factor $\Phi : \mathbb{X}_\alpha \to \mathbb{X}_\beta$ with the additional property that if $X$ is a Poisson point process on $G$ with intensity $\alpha$, then $X - \Phi(X)$ is a Poisson point process on $G$ with intensity $\alpha - \beta$.*

The constructions up until this point build on the ones in [8]. But toward proving Theorem 1.1 we apply tools from [22]. We note the method described for proving Theorem 1.2 is not bijective; when generating a Poisson point process on a bounded Borel set, the unique empty process occurs with nonzero probability. Our path around this obstruction is an intermediate isomorphism from a Poisson system to a particular product of Poisson systems, which we state in Proposition 1.5 after further definitions.

Let $\alpha, \beta_1, \beta_2 > 0$ and $\psi : \mathbb{X}_\alpha \to \mathbb{X}_{\beta_1} \times \mathbb{X}_{\beta_2}$. Denote $\psi(\cdot) = (\psi(\cdot)_1, \psi(\cdot)_2)$. If $\psi \circ g = (g \circ \psi_1, g \circ \psi_2)$ on a set of $P_\alpha$ full-measure, and for a Poisson point process $X$ on $G$ with intensity $\alpha$, the outputs $\psi(X)_i$ are independent Poisson processes on $G$ with intensities $\beta_i$ for $i = 1, 2$, then we say $\psi$ is a *factor from* $\mathbb{X}_\alpha$ *to* $(\mathbb{X}_{\beta_1}, \mathbb{X}_{\beta_2})$. The map $\psi$ is *finitary* if each coordinate mapping is finitary.

Likewise, let $\zeta : \mathbb{X}_{\beta_1} \times \mathbb{X}_{\beta_2} \to \mathbb{X}_\alpha$. If $\zeta \circ g = g \circ \zeta$ and $(P_{\beta_1} \times P_{\beta_2}) \circ \zeta^{-1} = P_\alpha$ on a set of $P_{\beta_1} \times P_{\beta_2}$ full-measure, then $\zeta$ is a *factor from* $(\mathbb{X}_{\beta_1}, \mathbb{X}_{\beta_2})$ *to* $\mathbb{X}_\alpha$. Let $\mu_1, \mu_2, \nu_1, \nu_2 \in \mathbb{M}$. A *coding window* of $\zeta$ is a function $\omega : \mathbb{M} \times \mathbb{M} \to \mathbb{N} \cup \{\infty\}$, such that if

$$(\mu_1, \mu_2)|_{B(e, \omega(\mu_1, \mu_2))} = (\nu_1, \nu_2)|_{B(e, \omega(\mu_1, \mu_2))},$$

then $\zeta(\mu_1, \mu_2)|_{B(e,1)} = \zeta(\nu_1, \nu_2)|_{B(e,1)}$. A factor from $(\mathbb{X}_{\beta_1}, \mathbb{X}_{\beta_2})$ to $\mathbb{X}_\alpha$ is *finitary* if there exists a coding window of $\zeta$ that is finite $P_{\beta_1} \times P_{\beta_2}$ almost surely. A factor from $\mathbb{X}_\alpha$ to $(\mathbb{X}_{\beta_1}, \mathbb{X}_{\beta_2})$ is an (finitary) isomorphism if its inverse is a (finitary) factor from $(\mathbb{X}_{\beta_1}, \mathbb{X}_{\beta_2})$ to $\mathbb{X}_\alpha$.

In Proposition 1.5, we return to the usual requirements on $G$.

PROPOSITION 1.5. *Let $\alpha, \beta > 0$. There exists a finitary isomorphism from $\mathbb{X}_\alpha$ to $(\mathbb{X}_\alpha, \mathbb{X}_\beta)$.*

With Proposition 1.5, we modify the construction used in the proof of Theorem 1.2 in order to prove Theorem 1.1. Instead of using randomness from a source system to populate a single Poisson process on $G$, we populate a pair of Poisson processes on $G$, bijectively. We make sure members of the pair are independent so that we may send the permuted pair back through the inverse factor map. All these operations together constitute an isomorphism between Poisson systems. As a corollary to Proposition 1.5 and Theorem 1.1, we gain the next result.

COROLLARY 1.6. *Poisson systems over $G$ are finitarily isomorphic to products of Poisson systems over $G$.*

We prove Theorems 1.1, 1.2 and 1.3, Proposition 1.5 and Corollary 1.6 in Section 6.2 (mostly they follow from work done in previous sections for specific cases of $G$), and prove Theorem 1.4, on Poisson splitting, in Section 7.

The main difficulties of the paper are ensuring the collection and distribution of randomness performed by factor maps uphold independence properties, and securing bijectivity for the isomorphism. Tools for dealing with the latter port easily from the $\mathbb{R}^n$ setting and [22], while dealing with the former is not as direct and informs the restrictions on $G$. We point out some of the differences between our construction and that of [8] as they arise.

**2. Poisson systems over groups.** Although Poisson point processes are most commonly defined on $\mathbb{R}^n$, all that is needed for existence is a Borel space with measurable singletons and a nonatomic measure satisfying a weak finiteness condition; see [15], Chapter 2. Considering a Poisson point process as a measure-preserving dynamical system restricts the possibilities somewhat.

Recall $\mathbb{M}$ is the space of Borel simple point measures on $G$. An element of $\mathbb{M}$ is a collection of "turned on" points (points with simple point measure value 1, simple meaning a point may only be turned on once) scattered throughout $G$. The scattering obeys the law $P_\alpha$ for some $\alpha > 0$. For $g \in G$ to act as a measure-preserving transformation on $\mathbb{M}$, the group $G$ must have a left-invariant measure. Haar's theorem states any locally compact group has such a measure, and a later converse proved in stages by Weil and Mackey states for a left-invariant measure to exist on a group, the group must be locally compact ([5], Chapter 1). Hence, for a Poisson system over $G$ to be defined, we cannot weaken the locally compact condition on $G$.

Further, the intensity of the process must be a constant value. Poisson point processes with constant intensity are called *homogeneous*, or for reasons we will address in this section, *uniform*.

If Haar measure on $G$ is $\sigma$-finite, then it satisfies the necessary conditions to host a Poisson process. It may be possible to remove the $\sigma$-finite condition and still have examples of groups where Poisson systems are valid, but this situation is less clear than the locally compact one. We note if a locally compact group is not $\sigma$-compact, its Haar measure is not $\sigma$-finite ([7], Chapter 2).

2.1. *Construction specific requirements.* Our approach to prove Theorem 1.1 is by construction. In order to construct a factor from $\mathbb{X}_\alpha$ to $\mathbb{X}_\beta$, we create an output Poisson point process from a source process. To do so, we need plenty of independent randomness from the source, which we collect with a marker system. The construction must be $G$-equivariant, limiting how we can define markers. The definition is based on distances between and densities of clustered points in the source process. Thus, we require that $G$ have a left-invariant metric.

Recall Birkhoff and Kakutani's theorem that a group is metrizable if and only if it is first countable. Struble proved a locally compact group is first countable and $\sigma$-compact if and only if it is second countable [24]. We require $G$ to be locally compact and $\sigma$-compact in order for its Haar measure to be left-invariant and $\sigma$-finite. Since we need a left-invariant metric, we reasonably require $G$ to be locally compact and second countable (in other words, locally compact Polish). We often directly use the second countable property (which implies hereditarily Lindelöf) in proofs to find countable open covers of $G$ and of open subsets of $G$, replacements for the $\mathbb{Z}^n$ lattices used in [8].

The marker system in [8] depends on a weak connectedness property of $\mathbb{R}^n$, stated at the start of Section 3. A topological group $T$ is said to be *compactly generated* if there exists a compact set $K \subseteq T$ such that
$$T = \bigcup_{i \geq 1} (K \cup K^{-1} \cup \{e\})^i.$$
In order for a locally compact Polish group to satisfy the marker property, it suffices that it be compactly generated. We discuss some nuances in Sections 3 and 6, and remark that the compactly generated condition is related to the notion of connectedness; see [5], Proposition 4.B.8.

A locally compact Polish group may not be compactly generated, but it must contain a compactly generated open subgroup $H$ of countable index ([5], Corollary 2.C.6). Then we may apply a version of Theorem 1.1 for compactly generated groups (Theorem 5.3, if $H$ is noncompact) to Poisson systems over $H$, as well as the other copies of $H$ in $G$ given by left cosets. This is Stepin's method of coinduction: A $G$-action is induced from the existing $H$-actions on cosets. Bowen uses coinduction in [3] to prove there exists a factor from any Bernoulli shift to any other one for countable groups containing a non-Abelian free subgroup and includes a proof of Stepin's theorem from [23].

The proof applies to our setting in the case when $H$ is noncompact with minimal adjustments, which we make in Section 5.2. If $G$ has no noncompact, compactly generated open subgroup, we do not apply coinduction for technical reasons. We consider this case and discuss technicalities in Section 6.1. As an example of a locally compact Polish group that is not compactly generated but contains a noncompact, compactly generated open subgroup, consider $H \times \mathbb{R}$, where
$$H := \bigoplus_{i \in \mathbb{N}} \mathbb{Z}/2\mathbb{Z}.$$
On the other hand, the group $H \times S^1$, where $S^1$ is the circle group, is locally compact Polish, but has no noncompact, compactly generated open subgroup.

At the core of our construction is a liberal use of the Borel isomorphism theorem for measures, stated in Remark 2.1 as found in [12], Theorem 17.41. We need Borel subsets of $G$ to be standard Borel spaces (i.e., isomorphic to a Polish space in which measurable sets are Borel) in order to apply the isomorphism theorem to extract randomness and populate an output Poisson process. Such is the case when $G$ is Polish ([12], Corollary 13.4).

We say a probability measure $\kappa$ on a measurable space $(A, \mathcal{B}(A))$ is *continuous* if $\kappa(\{a\}) = 0$ for all $a \in A$. Note Haar measure on a discrete, locally compact group is a counting measure. Thus, $G$ being nondiscrete is integral to our construction. Certainly, Theorem 1.1 is false for countable sofic groups.

REMARK 2.1 (Borel isomorphism theorem for measures). Let $A$ be a standard Borel space and $\kappa$ a continuous probability measure on $(A, \mathcal{B}(A))$. Then there exists a Borel isomorphism $S : A \to [0, 1]$ with $\kappa \circ S^{-1}$ Lebesgue measure on $[0, 1]$.

Occasionally, we must integrate over the space $\mathbb{M}$ equipped with the smallest $\sigma$-algebra such that the projection maps

$$\pi_A : \mathbb{M} \to \mathbb{N} \cup \{\infty\}$$

sending $\mu \in \mathbb{M}$ to $\mu(A)$ are measurable for all $A \in \mathcal{B}(G)$. We denote this $\sigma$-algebra as $\mathcal{M}(G, \mathcal{B}(G))$. For this purpose, we want Poisson processes on $G$, that is, random variables taking values on $\mathbb{M}$, to act as well behaved as they do on $\mathbb{R}^n$ and $\mathbb{M}(\mathbb{R}^n, \mathcal{B}(\mathbb{R}^n))$. This is the case when $G$ is a complete separable metric space (equivalently, Polish; see [18], Chapter 8). Then $\mathbb{M}$ is such a space as well. Further, Haar measure on $G$ is $\sigma$-finite. In particular, Fubini's theorem applies, useful for proving independence properties.

For a homogeneous Poisson point process $X$ on a complete separable metric space, given a bounded Borel set $A$ of the space, the points of $X|_A$ are uniformly distributed throughout $A$ ([18], Theorem 1.2.1). We use this property many times throughout the paper and state its most useful form to us in Remark 2.2.

Throughout the paper, we fix $\lambda$ as Haar measure on $G$. A uniform random variable without reference to a set is assumed to be on the interval $[0, 1]$. We denote the Dirac measure with mass at $g \in G$ as $\delta(g)$.

REMARK 2.2 (Distribution on finite measure sets). Let $X$ be a Poisson point process on $G$ with intensity $\alpha$, and let $A \in \mathcal{B}(G)$ such that $\lambda(A)$ is finite. Let $Z$ be a Poisson random variable with mean $\alpha \cdot \lambda(A)$, and let $\{U_i\}_{i \in \mathbb{N}}$ be a sequence of independent random variables, each uniformly distributed on $A$ and independent of $Z$. Then

$$X|_A \stackrel{d}{=} \sum_{i=1}^{Z} \delta(U_i).$$

It remains to address the noncompact requirement on $G$. A Poisson point process of intensity $\alpha$ on a compact group $K$ is the empty process $\varnothing \in \mathbb{M}(K)$ with probability $e^{-\alpha \cdot \lambda(K)} > 0$. Under our construction, it is impossible to produce an output process from an empty source process. But as long as $G$ is noncompact, by definition we are guaranteed, almost surely, a Poisson point process on $G$ with countably infinitely many points.

2.2. *Nonunimodular case.* Suppose $G$ is locally compact and $\sigma$-compact (so Poisson systems are well defined). Benjy Weiss proposed the following isomorphism when $G$ is nonunimodular and the modular homomorphism $\Delta : G \to (0, \infty)$ is surjective. (The modular homomorphism measures the difference between left and right Haar measures on $G$; see [7], Chapter 2 for a definition.) Fix $\alpha > 0$, and let $X$ be a Poisson point process on $G$ with intensity $\alpha$. Let $g \in G$. The map sending $X$ to $Xg$ is a finitary isomorphism from $\mathbb{X}_\alpha$ to $\mathbb{X}_{\alpha \cdot \Delta(g)}$. Since $\Delta$ is surjective, all Poisson systems over $G$ are isomorphic.

This gives an incredibly short proof of Theorem 1.1 under nonunimodular and surjective modular homomorphism assumptions.

2.3. *Extending the group action.* While we maintain the finitary property belonging to the factor maps in [8, 22] in our construction, we do not extend the group action from the base group. The referenced factor maps are between systems $(\mathbb{M}(\mathbb{R}^n), P_\alpha, \text{Isom}\,\mathbb{R}^n)$ and $(\mathbb{M}(\mathbb{R}^n), P_\beta, \text{Isom}\,\mathbb{R}^n)$, where $\text{Isom}\,\mathbb{R}^n$ is the isometry group of $\mathbb{R}^n$. Our factors are between $(\mathbb{M}(G), P_\alpha, G)$ and $(\mathbb{M}(G), P_\beta, G)$.

Beyond extending the group action to an isometry group, one may think of extensions more generally. For a locally compact group $G$ and a closed subgroup $H$ with modular homomorphisms $\Delta_G$ and $\Delta_H$, there exists a $G$-invariant Radon measure on $G/H$ if and only if $\Delta_G|H = \Delta_H$ ([7], Theorem 2.49). Lewis Bowen posed the following question.

QUESTION 1. Let $G$ be locally compact, let $H$ be a closed subgroup of $G$ such that $\Delta_G|H = \Delta_H$, and let $\alpha, \beta > 0$. Under what further conditions are $(\mathbb{M}(G/H), P_\alpha, G)$ and $(\mathbb{M}(G/H), P_\beta, G)$ isomorphic?

Evans shows in [6] that for the Poisson splitting constructed in [8] (related to but not precisely the finitary factor map of the same paper), it is impossible to extend equivariance further than $\text{Isom}\,\mathbb{R}^n$ (see also Question 3).

**3. Markers and Voronoi tessellations.** Our first step toward an isomorphism is to construct markers. To do so, we generalize the ones on $\mathbb{R}^n$ found in [8, 22]. These $\mathbb{R}^n$ markers rely on some connected-like behavior of $\mathbb{R}^n$ equipped with the Euclidean metric; namely, for all $g, h \in \mathbb{R}^n$ and $m \in \mathbb{N}$ with $|g - h| = m$ and $p, q \in \mathbb{N}$ such that $p + q = m$, there exists $k \in \mathbb{R}^n$ such that $|g - k| = p$ and $|h - k| = q$.

A word metric on any group always exhibits such behavior and is left-invariant, but depending on the group, may lack any number of desirable properties, such as measurability, properness or compatibility. However, a locally compact and compactly generated group has a measurable, proper word metric ([5], Remark 4.B.3) with other serviceable properties, which we identify as we use.

We fix $G$ as a nondiscrete, noncompact, locally compact and compactly generated Polish group throughout this section, Section 4 and Section 5.1, and $S$ as a compact generating subset of $G$ such that the topological boundary of any power of $S$ has measure zero (we define a topological boundary after Corollary 3.8, prove we may choose such $S$ in Lemma 3.9, and apply the requirement in Lemma 3.10).

Since $G$ is uncountable, so is $S$. The *word metric defined by $S$ on $G$* is the map $d_S : G \times G \to \mathbb{N}$ such that for any $g, h \in G$, we have

$$d_S(g, h) := \min\left\{ n \geq 0 : \begin{array}{l} \text{there exist } s_1, \ldots, s_n \in S \cup S^{-1} \\ \text{such that } g^{-1}h = s_1 \ldots s_n \end{array} \right\}.$$

We define the ball centered at $g \in G$ of radius $r \in \mathbb{N}$ as

$$B(g, r) := \{h \in G : d_S(g, h) \leq r\},$$

replacing our earlier definition. Note we include elements at distance $r$ from $g$. Balls are compact with respect to the underlying topology of $G$, since $B(e, m) = (S \cup S^{-1})^m$ for any $m \in \mathbb{N}$. The Borel set

$$A(g, a, b) := \{h \in G : a \leq d_S(g, h) \leq b\}$$

is the *shell centered at $g$ from $a$ to $b$*. Shells will be used to "mark" sources of randomness.

Before we continue with further definitions, we state a lemma important to our construction of a Poisson splitting. Lemma 3.1 appears as Lemma 12 in [8], and we refer to [8] for the proof.

LEMMA 3.1 (Poisson coupling). *Let $r \in (0, 1)$. There exists a constant $C_r$ such that for $\gamma > C_r$, there exist Poisson random variables $X$, $Y$ and $X + Y$ with means $r\gamma$, $(1 - r)\gamma$ and $\gamma$, respectively, and*

$$\mathbb{P}(Y = 0 | X + Y = 1) = \mathbb{P}(X = 0 | X + Y = 2) = 1.$$

The precise importance of Lemma 3.1 becomes clear during the splitting construction (in particular, in Proposition 7.3). For now, we let $\alpha > \beta > 0$, set $r := \beta/\alpha$ and fix $\rho \geq 10$ so that $\lambda(A(e, 6, \rho)) > C_r$ (there exists some $\rho$ large enough to satisfy the inequality; we require $\rho \geq 10$ specifically for enough space to prove what we need). We do not need this condition to prove Theorem 1.1, but it does not make the proof any more difficult.

Recall $N(A)$ is the number of points of $X$ in $A$, for any $A \in \mathcal{B}(G)$ with finite measure. We say $g \in G$ is a *seed* if $B(g, 3\rho + 30)$ satisfies the following properties:

1. For every ball $B$ of radius 1 such that $B \subset A(g, \rho + 6, \rho + 10)$, $N(B) \geq 1$.
2. The shell $A(g, \rho + 10, 3\rho + 30)$ is empty; that is, $N(A(g, \rho + 10, 3\rho + 30)) = 0$.

For any seed $g$, we call $A(g, \rho + 6, \rho + 10)$ its *dense shell* and $A(g, \rho + 10, 3\rho + 30)$ its *empty shell*. Seeds are defined so that we may place an equivalence relation on them, and those within a small enough distance of each other belong to the same equivalence class; the remaining ones must be far enough apart to avoid an intersection of dense and empty shells. In particular, we have chosen radii so seeds cannot spawn inside empty shells of other seeds, since $3\rho + 30 - (\rho + 10) > 2 \cdot (\rho + 6)$, and so we have a sizable enough gap between $g$ and $A(g, \rho + 6, \rho + 10)$ to prove independence and coupling properties. The particular constants are unimportant beyond meeting these requirements.

Our isomorphism construction deeply depends upon Lemma 3.2—it is either referenced directly or lurks in the background of most subsequent proofs.

LEMMA 3.2 (Seed distances). *Let $X$ be a Poisson point process on $G$ with intensity $\alpha$, and let $g, h \in G$ be seeds of $X$. Then $d_S(g, h) \notin [2, 4\rho + 38]$.*

PROOF. Note balls of radius 1 centered at elements in $A(g, \rho + 7, \rho + 9)$ contain at least one point of $X$, and $A(h, \rho + 10, 3\rho + 30)$ contains no points of $X$. We split the proof into two cases.

First, let $m, n \in \mathbb{N}$ such that $\rho + 7 \leq m \leq \rho + 9$ and $\rho + 11 \leq n \leq 3\rho + 29$, and suppose $d_S(g, h) = m + n$. Let $s_1, \ldots, s_{m+n} \in S \cup S^{-1}$ such that $g^{-1}h = s_1 \ldots s_{m+n}$. Then $d_S(g, gs_1 \ldots s_i) = i$ for $1 \leq i \leq m + n$. Thus, we may suppose $g' \in B(gs_1 \ldots s_m, 1)$ with $\delta(g') = 1$ (recall $\delta(g')$ is the Dirac measure with mass at $g'$). Consider that $h^{-1}gs_1 \ldots s_m = (s_{m+1} \ldots s_{m+n})^{-1}$, and

$$d_S(h, gs_1 \ldots s_m) \geq d_S(g, h) - d_S(g, gs_1 \ldots s_m) = n,$$

so $d_S(h, gs_1 \ldots s_m) = n$ by the triangle inequality. But then $g' \in A(h, \rho + 10, 3\rho + 30)$. We have a contradiction via a "turned on" point in the dense shell of $g$ belonging also to the empty shell of $h$, so $d_S(g, h) \notin [2\rho + 18, 4\rho + 38]$.

For the next case, fix $m = \rho + 9$, and let $n \in \mathbb{N}$ such that $\rho + 11 \leq n \leq 3\rho + 29$. Suppose $d_S(g, h) = n - m$. Let $s_1, \ldots, s_{n-m} \in S \cup S^{-1}$ such that $g^{-1}h = s_1 \ldots s_{n-m}$, and let $s'_1 \ldots s'_m \in S \cup S^{-1}$ such that $gs_1 \ldots s_{n-m}s'_1 \ldots s'_m =: h'$ is irreducible. Then $h^{-1}h' = s'_1 \ldots s'_m$ and $d(h, h') = m$. But $g^{-1}h' = s_1 \ldots s_{n-m}s'_1 \ldots s'_m$ and $d(g, h') = n$, so a turned on point in the dense shell of $h$ belongs to the empty shell of $g$, another contradiction. Thus, $d_S(g, h) \notin [2, 2\rho + 20]$. Combined with the first case, we are done. □

In light of Lemma 3.2, we define an equivalence relation $\sim$ on seeds by setting $g \sim h$ whenever $d_S(g, h) < 2$. Given a seed $g$, we denote its equivalence class of seeds under $\sim$ as $[g]$. The *core* of a seed class is the uncountable, closed set

$$\mathcal{C}[g] := \bigcap_{h \in [g]} B(h, 3).$$

Note cores have positive measure by Lemma 3.2, since $B(h, 1) \subseteq \mathcal{C}[g]$ for any $h \in [g]$. Some cores contain a unique point of $X$; we call such cores *identifiable*. The unique $X$-point in an identifiable core is its *landmark* $\ell[g] := X|_{\mathcal{C}[g]}$. We denote the union or sum (depending on our perspective of a landmark as a point in $G$ or a random variable, which we switch between) of landmarks as $\ell(X)$.

For any subset $A$ of $G$, set

$$\operatorname{diam} A := \sup_{g, h \in A} d_S(g, h).$$

We say $A$ is *bounded* if $\operatorname{diam} A < \infty$. Since $G$ is compactly generated, the boundedness of $A$ is independent of a choice between the metric $d_S$ and a geodesically adapted compatible metric ([5], Corollary 4.B.11). Remark 3.3 summarizes information from Definition 3.B.1, Corollary 4.B.11, Remarks 4.A.3 and 4.B.3, and Proposition 4.B.9 in [5].

REMARK 3.3 (Quasi-isometric metrics). On a locally compact and compactly generated group $G$, any two geodesically adapted metrics are quasi-isometric. In particular, the word metric $d_S$ is quasi-isometric to any left-invariant, proper, compatible metric on $G$, which implies any $d_S$-finitary factor is finitary.

We have $[g] \subseteq \mathcal{C}[g]$ and $\operatorname{diam} \mathcal{C}[g] \leq 3$ for any seed $g$. This informs our next definition. The *fitted shell* of a landmark is the uncountable Polish set

$$F[g] := A(\ell[g], 6, \rho),$$

and $F(X)$ refers to the union of fitted shells under $X$. Fitted shells are disjoint by Lemma 3.2. Each $X|_{F[g]}$ is a potential source of randomness. Lemma 3.4 helps to ensure independence conservation when collecting and distributing randomness from $X|_{F(X)}$.

LEMMA 3.4 (Fitted shell independence property). *Let* $\mu \in \mathbb{M}$, *and define*

$$\hat{F}(\mu) := \bigcup_{\ell[g] \in \ell(\mu)} A(\ell[g], 5, \rho + 1).$$

*Suppose* $\mu' \in \mathbb{M}$ *such that* $\mu|_{\hat{F}(\mu)^c} = \mu'|_{\hat{F}(\mu)^c}$. *Then* $\ell(\mu) = \ell(\mu')$ *and* $F(\mu) = F(\mu')$.

This property states that given two Poisson processes, which agree everywhere except on the fitted shells of one of the processes and relatively small "bumper shells" around them, both processes have the same fitted shells. In other words, the set $F(X)$ depends only on $X|_{\hat{F}(X)^c}$. Indeed, if two processes agree just on their cores, dense shells and empty shells, both processes should have the same fitted shells. We use the tolerance provided by the bumper shells in the proof of Lemma 4.2.

PROOF OF LEMMA 3.4. Let $g \in G$ be a seed under $\mu$ with identifiable core. Consider that $\mathcal{C}[g] \subseteq \hat{F}(\mu)^c$ and $\mu|_{\hat{F}(\mu)^c} = \mu'|_{\hat{F}(\mu)^c}$, so $\ell[g]$ is the landmark of $[g]$ under both $\mu$ and $\mu'$, as long as $g$ is a seed under $\mu'$. Thus, in order to show $F(\mu) = F(\mu')$, it suffices to restrict our attention to dense and empty shells and show

$$\mu|_{A(g, \rho+10, 3\rho+30)} = \mu'|_{A(g, \rho+10, 3\rho+30)}$$

to imply $g$ is a seed under $\mu'$. By Lemma 3.2, we know if $h$ is an identifiable seed under $\mu$ such that $h \notin [g]$, then $A(\ell[h], 5, \rho + 1) \cap B(g, 3\rho + 30) = \varnothing$. This, along with $\mu|_{\hat{F}(\mu)^c} = \mu'|_{\hat{F}(\mu)^c}$, gives us that $g$ is a seed under $\mu'$. □

Only those $X|_{F[g]}$ such that $N(X|_{F[g]}) \geq 1$ are isomorphic to uniform random variables. If $N(X|_{F[g]}) = 1$, then $X|_{F[g]}$ is uniformly distributed on $F[g]$, and the Borel isomorphism theorem from Remark 2.1 applies. If $N(X|_{F[g]}) > 1$, we could have an isomorphism still by the same theorem. But when $N(X|_{F[g]}) = 0$, the process $X|_{F[g]}$ is empty, and we cannot extract any useful randomness. In the language of the isomorphism theorem, the probability measure is not continuous.

To construct an output Poisson point process on $G$ given the input process $X$, we construct Poisson point processes on cells of a partition of $G$. The new processes as well as the partition both need to come from randomness within $X$, but the information sources should be independent of each other in order for the new processes to coalesce into a process on $G$. Thus, we want to place a partition on $G$ based on information independent of $X$ restricted to those fitted shells containing precisely one point of $X$ (we could use all fitted shells containing at least one point of $X$, but this approach streamlines construction slightly). To do so, we use a type of Voronoi tessellation.

The statement in Remark 3.5 is a corollary of the superposition theorem in [15] and justifies piecing processes together. Remark 3.5 holds for any complete separable metric space.

REMARK 3.5 (Poisson point process cut and paste). Let $X$ be a Poisson point process on $G$ with intensity $\alpha$, $N \subseteq \mathbb{N}$, and $\{A_i\}_{i \in N}$ a measurable partition of $G$. Let $\{Y_i\}_{i \in N}$ be a sequence of independent Poisson point processes, where each $Y_i$ is on $A_i$, with intensity $\alpha$. Then

$$\sum_{i \in N} Y_i \stackrel{d}{=} X.$$

We introduce the Voronoi partition for the factor and isomorphism constructions found in Section 5. Let $\mu \in \mathbb{M}$ and $x$ be a point of $\mu$. We use $x$ to express both an element of $G$ and a point measure. The *Voronoi cell* of $x$ is

$$\mathcal{V}(x) := \{g \in G : d_S(x, g) \leq d_S(y, g) \text{ for } y \in \mu\}.$$

In this context, we refer to each $x \in \mu$ as a *site*. The *Voronoi tessellation* of $\mu$ is the collection $\mathcal{V}(\mu) := \{\mathcal{V}(x)\}_{x \in \mu}$. Some elements of $G$ could belong to multiple Voronoi cells, although this cannot happen in $\mathbb{R}^n$ outside of measure zero sets. Abért and Mellick identify this distinction in [1] and introduce a tie-breaking function in pursuit of a true partition; we follow a similar approach. First, we define the particular points of a process we use as sites.

Given a seed $g$ with identifiable core, if its fitted shell contains precisely one process point, we say $g$ and its associated class, fitted shell, landmark and core are *harvestable*, with *harvest* $\mu|_{F[g]}$. The union of harvestable fitted shells is $F^*(\mu)$, and the union (or sum) of harvestable landmarks is $\ell^*(\mu)$. Any Poisson point process has "enough" harvestable landmarks for our construction $P_\alpha$ almost surely, which we prove in Lemma 3.6 with a Borel–Cantelli argument.

We tie-break cells with harvestable landmark sites in order to define our desired, measurable $G$-partition. We want to assign a unique value in $[0, 1]$ (almost surely) to each cell, to obtain a total ordering on cells. To make sure the tie-breaking map that does so is Borel, we define it in the same way regardless of the choice of cell. Some machinery is required.

Set $\mathcal{A} := \{(A, g) : A \text{ is closed in } B(e, 3) \subset G \text{ and } g \in A\}$. Note $(B(e, 3), d_S|_{B(e,3)})$ is a metric space and has a Hausdorff distance $d_H$ on its nonempty subsets; we recall its definition. For any nonempty $A \subseteq G$, the distance from $g \in G$ to $A$ is

$$d_S(g, A) := \inf_{a \in A} d_S(g, a).$$

For any nonempty $A, B \subseteq B(e, 3)$, the Hausdorff distance $d_H$ is defined as

$$d_H(A, B) := \max\left\{\sup_{a \in A} d_S(a, B), \sup_{b \in B} d_S(b, A)\right\}.$$

Now we define $d' : \mathcal{A} \times \mathcal{A} \to \mathbb{N}$ so that for $(A, g), (B, h) \in \mathcal{A}$,

$$d'((A, g), (B, h)) := d_H(A, B) + d_S(g, h).$$

The map $d'$ is a metric on $\mathcal{A}$. Define an equivalence relation $\sim$ on $\mathcal{A}$ by setting $(A, g) \sim (B, h)$ whenever there exists $f \in G$ such that $(A, g) = (fB, fh)$. The relation $\sim$ is closed in $\mathcal{A} \times \mathcal{A}$. It follows that the quotient space $\mathcal{A}/\sim$ with the quotient topology is compact and metrizable (since $(\mathcal{A}, d')$ is); in particular, the quotient space $\mathcal{A}/\sim$ is standard Borel. Finally, we fix a Borel isomorphism $T : \mathcal{A}/\sim \to [0, 1]$, and we set

$$T_{\mathcal{C}[g], \ell[g]} := T(\ell[g]^{-1}\mathcal{C}[g], e)$$

for all $\ell[g] \in \ell^*(\mu)$. The map $T$ is measurable, $G$-invariant and only depends on $\hat{F}(\mu)^c$ (as defined in Lemma 3.4). We use $T$ to define a modified Voronoi tessellation, which provides a true partition on $G$.

The *restricted Voronoi cell* of $\ell[g] \in \ell^*(\mu)$ is

$$\mathcal{V}^*(\ell[g]) := \left\{ f \in G : \begin{array}{l} \text{for each } \ell[h] \in \ell^*(\mu), \text{ either } d_S(\ell[g], f) < d_S(\ell[h], f) \text{ or} \\ d_S(\ell[g], f) = d_S(\ell[h], f) \text{ and } T_{\mathcal{C}[g], \ell[g]} \leq T_{\mathcal{C}[h], \ell[h]} \end{array} \right\},$$

and the *restricted Voronoi tessellation* of $\mu$ is $\mathcal{V}^*(\mu) := \{\mathcal{V}^*(\ell[g])\}_{\ell[g] \in \ell^*(\mu)}$. Restricted Voronoi tessellations are $G$-equivariant measurable partitions of $G$, almost surely.

In Lemma 3.6 and elsewhere, we use open covers made up of sets related to $d_S$-balls. The balls themselves are compact, but not necessarily open with respect to the underlying topology of $G$. If, for any $r \in \mathbb{N}$, we take the interior of $B(e, r + 1)$, we end up with an open set, which contains $B(e, r)$ (this can be shown by switching to a left-invariant, proper, compatible metric on $G$) and so have some constraint on the size of the interior.

For any $A \in \mathcal{B}(G)$, we denote the interior of $A$ as int $A$ and the closure of $A$ as cl $A$, both with respect to the topology of $G$.

LEMMA 3.6 (Infinitely many landmarks). *Let $X$ be a Poisson point process on $G$ with intensity $\alpha$, and let $A \in \mathcal{B}(G)$ such that $\lambda(A) = \infty$. Then $\mathbb{P}_\alpha$ almost surely $A$ contains infinitely many harvestable landmarks of $X$.*

PROOF. Fix $r \in \mathbb{N}$ such that $r > 6\rho + 60$. The collection $\{\text{int } B(g, r)\}_{g \in A}$ is an open cover of $A$. It is also an open cover of the open set

$$\bigcup_{g \in A} \text{int } B(g, r).$$

Since $G$ is second countable, it is also hereditarily Lindelöf, and we may let the collection $\{\text{int } B(g_i, r)\}_{i \in \mathbb{N}}$ be a countable subcover of $\{\text{int } B(g, r)\}_{g \in A}$ for some $\{g_i\}_{i \in \mathbb{N}} \subseteq A$.

Let $\{\text{int } B(g_{i_j}, r)\}_{j \in \mathbb{N}}$ be the disjoint subcollection of $\{\text{int } B(g_i, r)\}_{i \in \mathbb{N}}$ contained in $A$, where $g_{i_1} := g_k$ and $k$ is the minimal value such that int $B(g_k, r) \subseteq A$, $g_{i_2} := g_\ell$ and $\ell$ is the minimal value such that int $B(g_\ell, r) \subseteq A$ and int $B(g_k, r) \cap$ int $B(g_\ell, r) = \varnothing$, and the remaining $g_{i_j}$ are chosen similarly. Let $E_j$ be the event that int $B(g_{i_j}, r/2)$ contains a harvestable landmark of $X$.

It follows from the definition of a harvestable landmark that $\mathbb{P}(E_j) > 0$. By Lemma 3.2, the events $\{E_j\}_{j \in \mathbb{N}}$ are pairwise disjoint (the restriction on $r$ forces a distance of at least $6\rho + 60$ between any int $B(g_{i_j}, r/2)$). Furthermore, the $\mathbb{P}(E_j)$ are equal for all $j$. Then Borel–Cantelli implies the events $E_j$ occur infinitely often $P_\alpha$ almost surely. □

Later, we use the randomness from each harvest to populate a Poisson point process inside the restricted Voronoi cell of its landmark site. The following lemma and corollary are key in showing our constructions are finitary. For $g \in G$, we denote the unique harvestable landmark in the same restricted Voronoi cell as $g$ under $X$ as $\ell(X, g)$, and the restricted Voronoi cell itself as $\mathcal{V}(X, g)$.

LEMMA 3.7 (Restricted Voronoi cells have finite expected measure). *Let $X$ be a Poisson point process on $G$ with intensity $\alpha$ and restricted Voronoi tessellation $\mathcal{V}^*(X)$. Fix $\varepsilon > 0$ so that $d_S(\ell[g], \ell[h]) \geq 2\varepsilon$ for all $\ell[g], \ell[h] \in \ell^*(X)$ (we may do so by Lemma 3.2). Then*

$$\mathbb{E}\big(\lambda(\mathcal{V}(X, e))\mathbb{1}_{e \in B(\ell(X,e),\varepsilon)}\big)$$

*is finite $P_\alpha$ almost surely.*

PROOF. Let $\mu \in \mathbb{M}$ and $g \in G$. Define $M : \mathbb{M} \times G \to \mathbb{M} \times \mathbb{M}$ so that $M((\mu, g)) := (\mu, g\mu)$. Let $\nu \in \mathbb{M}$. Define $\Lambda : \mathbb{M} \times \mathbb{M} \to \{0, 1\}$ so that

$$\Lambda(\mu, \nu) := \begin{cases} 1 & \text{if } \nu = g\mu \text{ for some } g \text{ such that } g^{-1} \in B(\ell(\mu, e), \varepsilon), \\ 0 & \text{otherwise.} \end{cases}$$

Note $\Lambda(M(\mu, g)) = 1$ if and only if $g^{-1} \in B(\ell(\mu, e), \varepsilon)$. Thus, we have

$$\int_G \Lambda(M(\mu, g))\, d\lambda(g) = \lambda(B(\ell(\mu, e), \varepsilon)) = \lambda(B(e, \varepsilon)).$$

Further, integrating the above over $\mathbb{M}$ with respect to the probability measure $P_\alpha$ yields the same result. Since $\lambda$ is $\sigma$-finite and left-invariant, we have

$$\int_\mathbb{M} \int_G \Lambda(M(\mu, g))\, d\lambda(g)\, dP_\alpha(\mu) = \int_G \int_\mathbb{M} \Lambda(M(g^{-1}\mu, g))\, dP_\alpha(g^{-1}\mu)\, d\lambda(g)$$

$$= \int_\mathbb{M} \int_G \Lambda(M(g^{-1}\mu, g))\, d\lambda(g)\, dP_\alpha(\mu).$$

We claim the last integral is

$$\mathbb{E}\big(\lambda(\mathcal{V}(\mu, e))\mathbb{1}_{e \in B(\ell(\mu,e),\varepsilon)}\big).$$

To prove the claim, suppose $\Lambda(M(g^{-1}\mu, g)) = 1$. Then $g^{-1} \in B(\ell(g^{-1}\mu, e), \varepsilon)$, or equivalently, $e \in B(\ell(\mu, g), \varepsilon) \subseteq \mathcal{V}(\mu, g)$. This implies $g \in \mathcal{V}(\mu, e)$. On the other hand, if we suppose $e \in B(\ell(\mu, e), \varepsilon)$ and $g \in \mathcal{V}(\mu, e)$, then we have $g^{-1} \in B(\ell(g^{-1}\mu, e), \varepsilon)$, implying $\Lambda(M(g^{-1}\mu, g)) = 1$. So,

$$\int_G \Lambda(M(g^{-1}\mu, g))\, d\lambda(g) = \lambda(\mathcal{V}(\mu, e))\mathbb{1}_{e \in B(\ell(\mu,e),\varepsilon)}.$$

Integrating the above over $\mathbb{M}$ with respect to $P_\alpha$ proves the claim.

We have shown

$$\mathbb{E}\big(\lambda(\mathcal{V}(\mu, e))\mathbb{1}_{e \in B(\ell(\mu,e),\varepsilon)}\big) = \lambda(B(e, \varepsilon)) < \infty,$$

which proves the lemma. □

COROLLARY 3.8 (Bounded Voronoi cells). *Let $X$ be a Poisson point process on $G$ with intensity $\alpha$ and restricted Voronoi tessellation $\mathcal{V}^*(X)$. Each restricted Voronoi cell is bounded $P_\alpha$ almost surely.*

PROOF. Suppose for a contradiction there exists an unbounded restricted Voronoi cell of $X$ with positive probability under $P_\alpha$. Then the probability under $P_\alpha$ that $V := \mathcal{V}^*(\ell[g])$ is unbounded is positive, for any $\ell[g] \in \ell^*(X)$.

If $V$ is unbounded, for all $M \in \mathbb{N}$, there exists $f_M \in V$ such that $M \leq d_S(\ell[g], f_M) \leq d_S(\ell[h], f_M)$ for all $\ell[h] \in \ell^*(X)$. Note $B(f_M, M/2) \subseteq V$, so $\lambda(B(f_M, M/2)) \leq \lambda(V)$, implying $\lambda(V) = \infty$. But then

$$\mathbb{E}\big(\lambda(\mathcal{V}(X, e))\mathbb{1}_{e \in B(\ell(X,e),\varepsilon)}\big) = \infty,$$

contradicting Lemma 3.7. □

The *topological boundary* of a set $A \subseteq G$ is

$$\partial A := \operatorname{cl} A \setminus \operatorname{int} A.$$

We end the section by proving restricted Voronoi cells are close enough to compact for our needs, meaning the topological boundaries of restricted Voronoi cells almost surely have measure zero under $\lambda$, in Lemma 3.10. We use this property to measurably populate output Poisson point processes on restricted Voronoi cells. Recall we chose the compact generating set $S$ that defines our word metric $d_S$ to have topological boundary with measure zero. Lemma 3.9 justifies this choice.

LEMMA 3.9 (Generating set with measure zero topological boundary). *There exists a compact generating set $S \subset G$ such that $\lambda(\partial(S^n)) = 0$ for all $n \in \mathbb{N} \setminus \{0\}$.*

PROOF. The word metric $d_S$ is not compatible with the Borel topology on $G$. Let $d$ be a left-invariant, proper, compatible metric on $G$. Denote the closed ball centered at $g \in G$ of radius $r \in \mathbb{R}_{\geq 0}$ with respect to $d$ as $B_d(g, r)$, and denote the corresponding sphere as $S_d(g, r)$. Let $K$ be a compact generating set of $G$, let $t \in \mathbb{R}_{\geq 0}$ such that $K \subseteq B_d(e, t)$, let $t' > t$ and let $n \in \mathbb{N} \setminus \{0\}$.

Consider the set $\mathcal{R}_n := \{B_d(e, r)^n : t < r < t'\}$. If $\mathcal{R}_n$ is uncountable, all but countably many of $\partial(B_d(e, r)^n) \subseteq (\partial B_d(e, r))^n = S_d(e, r)^n$ must have measure zero, since

$$\lambda\left(\bigcup_{t < r < t'} S_d(e, r)^n\right) \leq \lambda\big(B_d(e, t')^n\big) < \infty.$$

If $\mathcal{R}_n$ is uncountable for all $n \in \mathbb{N} \setminus \{0\}$, then we may set $S := B_d(e, r)$ for some $r$ such that $\lambda(\partial(S^n)) = 0$ for all $n$.

If there exists $n$ such that $\mathcal{R}_n$ is countable, then there exist $r_1, r_2 \in \mathbb{R}_{\geq 0}$ such that $t < r_1 < r_2 < t'$ and $B := B_d(e, r_1)^n = B_d(e, r_2)^n$ for such $n$. But then $\operatorname{cl} B \subseteq B$ and $B \subseteq \operatorname{int} B$, so $\partial B = \varnothing$. In this case, set $S = B$. □

For convenience, we enumerate the set of harvestable cores under $X$ as $\{\mathcal{C}_i\}_{i \in \mathbb{N}}$, and the sets of corresponding harvestable landmarks, fitted shells, tie-breaking values and restricted Voronoi cells as $\{\ell_i\}_{i \in \mathbb{N}}$, $\{F_i\}_{i \in \mathbb{N}}$, $\{T_i\}_{i \in \mathbb{N}}$ and $\{\mathcal{V}_i\}_{i \in \mathbb{N}}$, respectively.

LEMMA 3.10 (Almost compact Voronoi cells). *Let $X$ be a Poisson point process on $G$ with intensity $\alpha$ and restricted Voronoi tessellation $\mathcal{V}^*(X)$. The topological boundary of each restricted Voronoi cell has measure zero $P_\alpha$ almost surely.*

PROOF. Let $d$ and $B_d$ be defined as in Lemma 3.9. We assume $S$ is a compact generating set for $G$ such that $\lambda(\partial(S^n)) = 0$ for all $n \in \mathbb{N} \setminus \{0\}$ and consider the topological boundary of a restricted Voronoi cell $\mathcal{V}_i$ for some $i \in \mathbb{N}$. Let $v \in \partial \mathcal{V}_i$. We claim $v$ lives in the topological boundary of a translation of $S^m$ for some $m \in \mathbb{N}$. To prove the claim, it suffices to show for every $\varepsilon > 0$ that there exist $x, y \in B_d(v, \varepsilon)$ with $x \notin \ell_i S^m$ and $y \in \ell_i S^m$.

So, let $\varepsilon > 0$. Because $v \in \partial \mathcal{V}_i$, we may fix $x, y \in B_d(v, \varepsilon)$ such that $x \notin \mathcal{V}_i$ and $y \in \mathcal{V}_i$. The proof of the claim splits into cases. First, suppose $d_S(y, \ell_i) < d_S(x, \ell_i)$. Then $y \in \ell_i S^m$ and $x \notin \ell_i S^m$ for some $m \in \mathbb{N}$ by the definition of $d_S$; thus, $v \in \partial \ell_i S^m$.

Next suppose $d_S(y, \ell_i) = d_S(x, \ell_i) = n$. Since $x \notin \mathcal{V}_i$, there exists $j \in \mathbb{N}$ such that $x \in \mathcal{V}_j$. Almost surely, either $T_i < T_j$ or vice versa. First, assume $T_i < T_j$. Then $d_S(x, \ell_j) < n$ and $d_S(y, \ell_j) \geq n$, so $x \in \ell_j S^{n-1}$ and $y \notin \ell_j S^{n-1}$, implying $v \in \partial \ell_j S^{n-1}$. If $T_i > T_j$, then $d_S(x, \ell_j) \leq n$ and $d_S(y, \ell_j) > n$, so $x \in \ell_j S^n$, $y \notin \ell_j S^n$ and $v_i \in \partial \ell_j S^n$.

We cannot have $d_S(y, \ell_i) > d_S(x, \ell_i)$, so we have shown

$$\partial \mathcal{V}_i \subseteq \bigcup_{j,m \in \mathbb{N}} \partial \ell_j S^m,$$

where each $\partial \ell_j S^m$ has measure zero. Hence, $\lambda(\partial \mathcal{V}_i) = 0$ for all $i \in \mathbb{N}$, $P_\alpha$ almost surely. □

**4. Harvesting randomness.** In this section, we explore the independence properties necessary to our constructions. For any Poisson point process $X$ on $G$ and $A \in \mathcal{B}(G)$, the restrictions $X|_A$ and $X|_{A^c}$ are independent Poisson processes. Further, their sum has the same law as $X$ (see Remark 3.5). Introducing a Poisson point process $Y$ on $G$, independent of $X$ and of the same intensity, allows us to rephrase this fact in the following way, more conducive to proving independent properties:

$$Y|_A + X|_{A^c} \stackrel{d}{=} X.$$

We want a version of this fact for a random Borel set in $G$, in particular the union of fitted shells $F(X)$, which depends on $X|_{\hat{F}(X)^c}$. Lemma 3.4 in Section 3 is a first step in this direction. We need Lemma 4.2 as well to prove our desired version, which we state in Proposition 4.1.

PROPOSITION 4.1 (Fitted shell cut and paste). *Let $X$ and $Y$ be independent Poisson point processes on $G$, both with intensity $\alpha$. The process $Z := Y|_{F(X)} + X|_{F(X)^c}$ is equal in distribution to $X$, and $F(Z) = F(X)$.*

Proposition 4.1 restates Proposition 16 in [8], and our proof adapts the proofs of Proposition 16 and Lemma 18 in the same paper from Poisson processes on $\mathbb{R}^n$ to Poisson processes on $G$. As a corollary (Corollary 4.3), we get a similar statement for the union of harvestable fitted shells $F^*(X)$, which is really what we are working toward. Ignoring technical details, the main tools in the proof of Proposition 4.1 are Lemma 3.4, which tells us $F(X)$ depends only on $X|_{\hat{F}(X)^c}$, and Fubini's theorem.

Lemma 4.2, adapted from Lemma 17 in [8], is used in the proof of Proposition 4.1 and addresses some of those technical details.

LEMMA 4.2. *Let $X$ be a Poisson point process on $G$ with intensity $\alpha$ and $A$ be a bounded Borel set. There exists a finite set $S$, a collection of bounded Borel sets $\{F_s\}_{s \in S}$ and a collection of disjoint events $\{E_s\}_{s \in S}$ such that*:

1. *The union over $S$ of the events $E_s$ occurs $P_\alpha$ almost surely.*
2. *For all $s \in S$, conditioned on $E_s$, we have $F(X) \cap A \subseteq F_s \subseteq \hat{F}(X)$.*
3. *Let $\sigma(X|_{F_s^c})$ be the $\sigma$-algebra generated by $X|_{F_s^c}$. For all $s \in S$, the event $E_s$ is $\sigma(X|_{F_s^c})$-measurable.*

It would be more straightforward to apply Lemma 4.2 in the proof of Proposition 4.1 if, instead of relying on bounded $A \in \mathcal{B}(G)$, we could replace $A$ directly with $G$. We cannot do so because of the first property. The cardinality of the set $S$ matches the cardinality of some cover of $A$. Requiring $A$ to be bounded, let us find a finite subcover. If $S$ were even countably infinite, the probability of the union over $S$ of events indexed by $S$, each with probability less than one, would be zero almost surely.

In the proof, we use the fact that since $d_S$ is a proper metric quasi-isometric to a proper, compatible metric on $G$, subsets of finite diameter in $(G, d_S)$ are relatively compact (see Remark 3.3).

PROOF OF LEMMA 4.2. We consider landmarks under $X$ that live in $A$; that is, landmarks belonging to the set $\ell(X) \cap A$. The collection $\{B(g, \rho)\}_{g \in A}$ contains the fitted shells of such landmarks. Let $\mathcal{O}$ be the union of sets in the collection. (The set of fitted shells of those landmarks in $\ell(X) \cap A$ may not be contained in $A$ itself.) Note $\operatorname{cl} \mathcal{O}$ is compact.

The collection $\{\operatorname{int} B(g, 1)\}_{g \in \operatorname{cl} \mathcal{O}}$ is an open cover of $\operatorname{cl} \mathcal{O}$; let $\{\operatorname{int} B(g_i, 1)\}_{i=1}^N$ be a finite subcover, and let $\{D_i\}_{i=1}^N$ be the disjointification of the finite subcover, meaning $D_1 := \operatorname{int} B(g_1, 1)$ and for $i > 1$,

$$D_i := \operatorname{int} B(g_i, 1) \setminus \bigcup_{j=1}^{i-1} \operatorname{int} B(g_j, 1).$$

Let $E_i$ be the event that $D_i$ contains a landmark in $\ell(X) \cap A$. Although the $D_i$ are disjoint, the events $E_i$ are not, due to the nature of landmarks. To obtain disjoint events, for all $s \in \{0, 1\}^N$, define

$$E_s := \left( \bigcap_{1 \leq i \leq N : s_i = 1} E_i \right) \cap \left( \bigcap_{1 \leq i \leq N : s_i = 0} E_i^c \right)$$

and set $S := \{s \in \{0, 1\}^N : \mathbb{P}(E_s) > 0\}$. Now for all $s \in S$, the events $E_s$ are disjoint. Furthermore, the $E_s$ satisfy (1).

To prove (2), we want to define $F_s$ based on $E_s$ so that $F(X) \cap A \subseteq F_s \subseteq \hat{F}(X)$. Suppose $\ell[h] \in \ell(X) \cap D_i$: Then $d(g_i, \ell[h]) \leq 1$. For any $f$ in the fitted shell of $\ell[h]$, we have $5 \leq d(g_i, f) \leq \rho + 1$, so set

$$F_s := \bigcup_{1 \leq i \leq N : s_i = 1} A(g_i, 5, \rho + 1),$$

which satisfies our conditions.

Property (3) assures us $E_s$ does not depend on any information from $X|_{F_s}$ and follows from Property (2) and Lemma 3.4. □

We are ready to prove Proposition 4.1. Let $\Omega$ be a set. For any subset $A \subseteq \Omega$, we use $\mathbb{1}_A$ to denote the indicator function of $A$.

PROOF OF PROPOSITION 4.1. By Remark 3.5, for any fixed $B \in \mathcal{B}(G)$, we know

$$Y|_B + X|_{B^c} \stackrel{d}{=} X.$$

Our goal is to prove the same for the random union of fitted shells $F(X)$. First note, since $Z|_{F(X)^c} = X|_{F(X)^c}$, it is immediate from Lemma 3.4 that $F(Z) = F(X)$.

Fix a bounded $A \in \mathcal{B}(G)$ and let $S$, $\{F_s\}$, and $\{E_s\}$ be defined as in Lemma 4.2, under $X$. Fix $s \in S$ and let $\mathcal{A}$ be any measurable set of point measures in $\mathcal{M}$. We claim

(∗) $\qquad \mathbb{P}(\{X|_{F_s \cap A} + X|_{F_s^c \cap A} \in \mathcal{A}\} \cap E_s)$

is equal to

(**) $$\mathbb{P}(\{Y|_{A \cap F(X)} + X|_{A \setminus F(X)} \in \mathcal{A}\} \cap E_s).$$

Note $(*) = \mathbb{P}(\{X|_A \in \mathcal{A}\} \cap E_s)$ and $(**) = \mathbb{P}(\{Z|_A \in \mathcal{A}\} \cap E_s)$.

Let $\mu$ be the law of $X|_{F_s \cap A}$ and $\nu$ the law of $(X|_{F_s^c \cap A}, \mathbb{1}_{E_s}, A \cap F(X))$. The process $X|_{F_s}$ is independent of $(X|_{F_s^c}, \mathbb{1}_{E_s}, F(X))$ by Lemmas 3.4 and 4.2. Hence, $X|_{F_s \cap A}$ is independent of $(X|_{F_s^c \cap A}, \mathbb{1}_{E_s}, A \cap F(X))$. Let $x, y \in (\mathbb{M}, \mathcal{M})$ under $\mu$ and $(w, e, b) \in (\mathbb{M}, \mathcal{M}) \times (\mathbb{M}, \mathcal{M}) \times (G, \mathcal{B}(G))$ under $\nu$. We have

$$(*) = \iint \mathbb{1}_{[x+w \in \mathcal{A}]} \mathbb{1}_{E_s} \, d\mu(x) \, d\nu(w, e, b)$$

$$= \iint \mathbb{1}_{[x|_b + x|_{b^c} + w \in \mathcal{A}]} \mathbb{1}_{E_s} \, d\mu(x) \, d\nu(w, e, b)$$

$$= \int \left( \int \mathbb{1}_{[x|_b + x|_{b^c} + w \in \mathcal{A}]} \, d\mu(x) \right) \mathbb{1}_{E_s} \, d\nu(w, e, b)$$

$$= \iiint \mathbb{1}_{[y|_b + x|_{b^c} + w \in \mathcal{A}]} \mathbb{1}_{E_s} \, d\mu(x) \, d\mu(y) \, d\nu(w, e, b) = (**)$$

by Remark 3.5 and Fubini's theorem, proving the claim. Taking the sum over $S$ on both sides of $(*) = (**)$ yields $\mathbb{P}(X|_A \in \mathcal{A}) = \mathbb{P}(Z|_A \in \mathcal{A})$ by Property (1) of Lemma 4.2.

Now, to prove

$$X \stackrel{d}{=} Z,$$

fix $r \in \mathbb{N} \setminus \{0\}$ and consider the open $G$-cover $\{\text{int } B(g, r)\}_{g \in G}$. Let $\{\text{int } B(g_i, r)\}_{i \in \mathbb{N}}$ be a countable subcover and $\{D_i\}_{i \in \mathbb{N}}$ its disjointification (as defined in Lemma 4.2). For any finite N, the set $D_0 \cup \cdots \cup D_N$ is bounded, and we have

$$\mathbb{P}(X|_{D_0 \cup \cdots \cup D_N} \in \mathcal{A}) = \mathbb{P}(Z|_{D_0 \cup \cdots \cup D_N} \in \mathcal{A}).$$

Taking the limit as $N \to \infty$ on both sides gives us that $X$ and $Z$ are equal in distribution. □

With Proposition 4.1, we have shown $X|_{F(X)}$ and $X|_{F(X)^c}$ are independent Poisson point processes, conditioned on $F(X)$. For the factor and isomorphism maps, we want to populate each restricted Voronoi cell $\mathcal{V}^*(\ell[g])$ with a Poisson point process via randomness from the harvest $X|_{F[g]}$. This requires more information than we can get from $F(X)$ alone; we also need to know which fitted shells contain exactly 1 point of $X$. That is, we need to know $F^*(X)$.

If we know $F^*(X)$, we still do not gain any information regarding locations of points in harvestable shells. We do gain information on the process restricted to $F(X) \setminus F^*(X)$, but this does not disrupt the independence properties we ultimately require, which we prove in Corollaries 4.3 and 4.4, extensions of Lemma 24 in [8] and Corollaries 9 and 10 in [22].

COROLLARY 4.3 (Harvestable fitted shell cut and paste). *Let $X$ and $Y$ be independent Poisson point processes on $G$, both with intensity $\alpha$. The process $Z := Y|_{F^*(X)} + X|_{F^*(X)^c}$ is equal in distribution to $X$, and $F^*(Z) = F^*(X)$.*

PROOF. Since $Z|_{F^*(X)^c} = X|_{F^*(X)^c}$, it must be true that $Z|_{F(X)^c} = X|_{F(X)^c}$, which implies $F(Z) = F(X)$ by Lemma 3.4. Then of course $F^*(Z) = F^*(X)$. It remains to show $X|_{F^*(X)}$ is independent of $X|_{F^*(X)^c}$. Note

$$F^*(X)^c = F(X)^c \cup (F(X) \setminus F^*(X)).$$

Enumerate the set of all fitted shells as $\{F'_i\}_{i\in\mathbb{N}}$. From Proposition 4.1, we know the random variables $\{N(F'_i)\}_{i\in\mathbb{N}}$ are independent of each other and of $X|_{F(X)^c}$. It follows that the processes $X|_{F(X)\setminus F^*(X)}$ and $X|_{F^*(X)}$ are independent. Thus, $X|_{F(X)^c}$ and $X|_{F^*(X)}$ are independent as well. The set $F(X) \setminus F^*(X)$ depends only on $(X|_{F(X)^c}, \{N(F'_i)\}_{i\in\mathbb{N}})$, so we are done. □

With Corollary 4.3, we prove Corollary 4.4, which we use to prove that given a Poisson point process $X$ with intensity $\alpha$, we may generate processes on cells of $\mathcal{V}^*(X)$, which sum to a Poisson process on $G$ with intensity $\beta$. Harvests of $X$ generate the new processes, while $\mathcal{V}^*(X)$ depends on $\ell^*(X) \subseteq X|_{F^*(X)^c}$ and tie-breaking values, both independent of $X|_{F^*(X)}$.

For $A \subseteq G$ and $V$ a $\lambda$-uniform random variable taking values in $A$, we refer to $V$ as a *uniform random variable on $A$*.

COROLLARY 4.4 (Source independence). *Let $X$ be a Poisson point process on $G$ with intensity $\alpha$, and conditioned on $F^*(X)$, let*

$$U := \{U_{F_i}\}_{F_i \in F^*(X)}$$

*be a sequence of independent uniform random variables, independent of $X$. There exists a measurable map $\xi : F^*(X) \to [0, 1]$ such that*

$$\left(X|_{F^*(X)^c}, F^*(X), \{\xi(F_i)\}_{F_i \in F^*(X)}\right) \stackrel{d}{=} \left(X|_{F^*(X)^c}, F^*(X), U\right).$$

PROOF. Each harvest $X|_{F_i}$ is a Poisson point process with intensity $\alpha$ on $F_i$ conditioned to contain precisely 1 point of $X$. By Corollary 4.3, elements of the set of all harvests of $X$ are independent as well as independent of $X|_{F^*(X)^c}$, and by Lemma 3.6, the cardinality of the set of all harvests is $P_\alpha$ almost surely countably infinite. Further,

$$\mathcal{S}(X) := \left\{(\ell_i)^{-1} X|_{F_i}\right\}_{i\in\mathbb{N}}$$

is a sequence of independent and identically distributed Poisson processes on the set $A := A(e, 6, \rho)$. Each process in $\mathcal{S}(X)$ forgets its starting location, so $\mathcal{S}(X)$ is independent of $(X|_{F^*(X)^c}, F^*(X))$.

Let $V$ be a uniform random variable on $A$, independent of $X$. Fix a Borel isomorphism $R : A \to [0, 1]$ so that $R(V)$ is a uniform random variable on $[0, 1]$. Denote the shifted process point of $(\ell_i)^{-1} X|_{F_i}$ as $v_i$. Define $\xi : F^*(X) \to [0, 1]$ by $\xi(F_i) := R(v_i)$ for $F_i \in F^*(X)$, so $\xi$ is measurable.

For each $F_i \in F^*(X)$, we have that $\xi(F_i)$ is uniformly distributed on $[0, 1]$ and the sequence of independent uniform random variables $\{\xi(F_i)\}_{F_i \in F^*(X)}$ is independent of $(X|_{F^*(X)^c}, F^*(X))$. Thus,

$$\left(X|_{F^*(X)^c}, F^*(X), \{\xi(F_i)\}_{F_i \in F^*(X)}\right) \stackrel{d}{=} \left(X|_{F^*(X)^c}, F^*(X), U\right)$$

as desired. □

We have proven all independence properties necessary to the proof of Theorem 1.2. One further independence property is required to prove Proposition 1.5 and Theorem 1.1. We need to know the Voronoi partition and the output system are independent. Lemma 4.5 is a version of Lemma 16 in [22], and in fact is implied from the same lemma, which relies on Fubini's theorem, along with Corollaries 4.3 and 4.4.

LEMMA 4.5 (Partition and output independence). *Let $(A, \mathcal{A})$ be a measure space, and let $X$ be a random variable taking values in $A$. Let $\{\mathcal{V}(X)_i\}_{i \in \mathbb{N}}$ be a random measurable partition of $G$ dependent on $X$, and let $V$ be a uniform random variable. Suppose $\kappa : [0, 1] \times A \times \mathbb{N} \to \mathbb{M}$ is a measurable map such that for all $a \in A$ and $i \in \mathbb{N}$, $\kappa(V, a, i)$ is a Poisson point process on $\mathcal{V}(X)_i$ with intensity $\beta > 0$. Let $\{U_i\}_{i \in \mathbb{N}}$ be a sequence of independent uniform random variables independent of $X$. Then*

$$\sum_{i \in \mathbb{N}} \kappa(U_i, X, i)$$

*is a Poisson point process on $G$ with intensity $\beta$, independent of $X$.*

**5. Compactly generated cases.** We start in Section 5.1 with a version of Theorem 1.2 for compactly generated groups. We then prove versions of Proposition 1.5 and Theorem 1.1. Once the compactly generated case is complete, we extend our results for this case to nondiscrete, noncompact, locally compact Polish groups with a noncompact, compactly generated open subgroup in Section 5.2 using coinduction (the group itself does not have to be compactly generated here).

5.1. *Compactly generated group.* To prove Theorem 1.2 for compactly generated $G$, we construct a factor from $\mathbb{X}_\alpha$ to $\mathbb{X}_\beta$ via the harvests and restricted Voronoi tessellation of the source system $\mathbb{X}_\alpha$.

We note the factor and isomorphism constructions in [22] are more explicit than those in this section. The Voronoi cells which partition $\mathbb{R}^n$ in [22] are convex, bounded polytopes and, given uniform random variables, may be populated with Poisson point processes without resorting to the Borel isomorphism theorem. Our restricted Voronoi cells on $G$ lack any such unifying structure. To populate our cells, we require a careful application of the Borel theorem which relies on Lemma 3.10. Let $C(G)$ be the set of compact subsets of $G$.

THEOREM 5.1 (Finitary factor for compactly generated groups). *Let $G$ be a nondiscrete, noncompact, locally compact and compactly generated Polish group and let $\alpha, \beta > 0$. There exists a finitary factor from $\mathbb{X}_\alpha$ to $\mathbb{X}_\beta$.*

PROOF. Suppose $X$ and $Y$ are Poisson point processes on $G$ with intensities $\alpha$ and $\beta$, respectively. We construct a factor map $\phi : \mathbb{X}_\alpha \to \mathbb{X}_\beta$ so that

$$\phi(X) \stackrel{d}{=} Y.$$

Recall the restricted Voronoi tessellation of $X$ is $\mathcal{V}^*(X) = \{\mathcal{V}_i\}_{i \in \mathbb{N}}$, a measurable $G$-partition $P_\alpha$ almost surely. Each cell of $\mathcal{V}^*(X)$ is $P_\alpha$ almost surely bounded (Corollary 3.8) with measure zero topological boundary (Lemma 3.10). We populate a Poisson point process with intensity $\beta$ on each shifted cell $\mathcal{S}_i := (\ell_i)^{-1} \mathcal{V}_i$ with randomness from $X|_{F_i}$. Cells are shifted to avoid disrupting $G$-equivariance of our eventual factor map.

The population map must be measurable. As with our definition for the Voronoi cell tie-breaking map, we need to ensure this map does not depend on the choice of $\mathcal{S}_i$. The Borel isomorphism theorem implies a uniform random variable is isomorphic to a nonempty Poisson process on $\mathcal{S}_i$, but we need a single map that works for all shifted cells. Hence, we consider the set

$$\mathcal{K} := \{(K, \mu) : K \text{ is compact in } G \text{ and } \mu \in \mathbb{M}(K, \{K \cap A : A \in \mathcal{B}(G)\})\}$$

(this set is a Polish space; see [12], Theorem 4.25 and [18], Chapter 8). Additionally, we need to split the behavior of the population map according to whether $\mathcal{S}_i$ is "populated" with an empty process (in this case, we lack an isomorphism to a uniform random variable). So, let

$$\mathcal{K}_0 := \{(K, \varnothing) : K \text{ is compact in } G\} \subseteq \mathcal{K}.$$

Let $K \in C(G)$ and $\mu \in \mathbb{M}(K)$. Fix a Borel isomorphism $\tau' : \mathcal{K} \setminus \mathcal{K}_0 \to (0, 1]$, and define $\tau : \mathcal{K} \to [0, 1]$ by

$$\tau((K, \mu)) := \begin{cases} \tau'(K, \mu) & \text{whenever } \mu \neq \varnothing, \\ 0 & \text{otherwise.} \end{cases}$$

Since $\tau$ is not explicit, we need to supply further structure so that a proper number of cells receive an empty process. Denote the law of a Poisson point process on $K$ with intensity $\beta$ as $P_K$ and the push-forward measure of $P_K$ under $\tau|_{\{K\} \times \mathbb{M}(K)}$ as $\nu_K$. Then

$$\nu_K(0) = P_K(N(K) = 0) =: p_{K,0}.$$

Let $x \in [0, 1]$, and define $r_K : [0, 1] \to [0, 1]$ by $r_K(x) := \nu_K([0, x])$. Note $r_K([0, x]) = [p_{K,0}, \nu_K([0, x])]$. Hence, $r_K$ pushes $\nu_K$ on $[0, 1]$ forward to the sum of Dirac mass on $p_{K,0}$ of measure $p_{K,0}$ and Lebesgue measure on $(p_{K,0}, 1]$. Define $\psi : \mathcal{K} \to C(G) \times [0, 1]$ by $\psi((K, \mu)) := (K, r_K \circ \tau)$ and $\pi : C(G) \times [0, 1] \to \mathcal{K}$ by

$$\pi((K, x)) := \begin{cases} \psi^{-1}(K, x) & \text{whenever } x > p_{K,0}, \\ (K, \varnothing) & \text{otherwise.} \end{cases}$$

The map $\pi$ is measurable. Its first output coordinate $\pi((K, x))_1$ is $K$, and its second output coordinate $\pi((K, x))_2$ is a Poisson point process of intensity $\beta$ on $K$, for $x$ chosen uniformly from $[0, 1]$.

Let $\xi : F^*(X) \to [0, 1]$ be defined as in Corollary 4.4. Recall $\xi$ extracts randomness from the shifted harvestable fitted shells $\{(\ell_i)^{-1} F_i\}_{i \in \mathbb{N}}$ of $X$. By Corollary 4.4, $\{\xi(F_i)\}_{F_i \in F^*(X)}$ is a sequence of independent uniform random variables independent of $X_{F^*(X)^c}$. Set $\pi_i := \pi((\operatorname{cl} \mathcal{S}_i, \xi(F_i)))_2$ for $F_i \in F^*(X)$. Since $\mathcal{S}_i$ is almost surely bounded (Corollary 3.8) and $d_S$ is proper, $\operatorname{cl} \mathcal{S}_i$ is almost surely compact, so $\pi_i$ is well defined. By Lemma 3.10, each $\pi_i$ is almost surely a Poisson process of intensity $\beta$ on $\operatorname{int} \mathcal{S}_i$, or equivalently $\mathcal{S}_i$, independent of $X_{F^*(X)^c}$.

Define $\phi : \mathbb{X}_\alpha \to \mathbb{X}_\beta$ as

$$\phi(X) := \sum_{F_i \in F^*(X)} \ell_i \pi_i.$$

The landmarks $\ell^*(X)$ and tessellation $\mathcal{V}^*(X)$ depend only on $X_{F^*(X)^c}$ by definition. Hence, Remark 3.5 implies $\phi(X)$ is a Poisson point process of intensity $\beta$ on $G$.

We check $\phi$ is $G$-equivariant. Let $g \in G$. If $g$ acts on $X$, then shifted fitted shells remain the same:

$$(g\ell_i)^{-1} g F_i = (\ell_i)^{-1} F_i.$$

Similarly, shifted restricted Voronoi cells $\mathcal{S}_i$ are $g$-invariant. Thus,

$$\phi \circ g(X) = \sum_{F_i \in F^*(X)} (g\ell_i) \pi_i = g \circ \phi(X).$$

Corollary 3.8 implies $\phi$ is $d_S$-finitary, since each Voronoi cell under $X$ is completely determined by $X$ restricted to some bounded set containing the cell and adjacent cells, and each harvestable fitted shell $F_i$ is contained in the Voronoi cell $\mathcal{V}_i$. Then Remark 3.3 implies $\phi$ is finitary. $\square$

The factor map in Theorem 5.1 fails to be bijective. The law of a Poisson point process on a bounded Borel set is not a continuous measure (the empty process is an atom), and so we cannot expect an isomorphism between such a process and a uniform random variable.

To construct a bijective factor map, we must preserve the information otherwise lost when a Poisson point process on a bounded Voronoi cell is determined to be empty, as well as those ambient points of the input process that do not determine the output process.

Our solution is the same as in [22] for the $\mathbb{R}^n$ case. Along with generating an output Poisson point process from an input one, we generate an additional output process using the "extra" information. We need to do some further work to end up with a final output of a single Poisson process, but this motivates Proposition 1.5. Proposition 5.2 is a version of Proposition 1.5 for our current requirements on $G$.

PROPOSITION 5.2. *Let $G$ be a nondiscrete, noncompact, locally compact and compactly generated Polish group and let $\alpha, \beta > 0$. There exists a finitary isomorphism from $\mathbb{X}_\alpha$ to $(\mathbb{X}_\alpha, \mathbb{X}_\beta)$.*

PROOF. Suppose $X$ and $Y$ are independent Poisson point processes on $G$ with intensities $\alpha$ and $\beta$, respectively. We build a factor map $\phi : \mathbb{X}_\alpha \to \mathbb{X}_\alpha \times \mathbb{X}_\beta$. We denote $\phi(X)$ as $(\phi(X)_1, \phi(X)_2)$ (and use similar notation for all coordinate maps) and construct $\phi$ so that

$$\phi(X) \stackrel{d}{=} (X, Y).$$

To do so, we modify the argument of Theorem 5.1. We populate a Poisson point process with intensity $\beta$ on each shifted cell $\mathcal{S}_i := (\ell_i)^{-1} \mathcal{V}_i$ with randomness from $X|_{F_i}$, but now we keep $X|_{\mathcal{S}_i}$ as well, except for on $F_i \subseteq \mathcal{S}_i$, where we resample the process.

Let $\mathcal{K}$, $\mathcal{K}_0$, $C(G)$, $\tau'$, $\tau$, $P_K$, $\nu_k$, $r_K$, $p_{K,0}$ and $\psi$ be defined as in Theorem 5.1, and let $K \in C(G)$. We need one more map before we adjust the definition of the population map. Let $x \in [0,1]$ and denote its binary expansion as $.x_1 x_2 x_3 \ldots$. Define $b: [0,1] \to [0,1]^2$ by $b(x) := (x^1, x^2)$, where $x^1 := .x_1 x_3 x_5 \ldots$ and $x^2 := .x_2 x_4 x_6 \ldots$. Let $U$ be a uniform random variable. The outputs $b(U)_1$ and $b(U)_2$ are independent uniform random variables.

Define $\pi : C(G) \times [0,1] \to \mathcal{K} \times [0,1]$ by

$$\pi((K,x)) := \begin{cases} \left(K, \varnothing, \dfrac{x}{p_{k,0}}\right) & \text{whenever } x < p_{K,0}, \\ \left(\psi^{-1}\left(K, b\left(\dfrac{x - p_{K,0}}{1 - p_{K,0}}\right)_1\right), b\left(\dfrac{x - p_{K,0}}{1 - p_{K,0}}\right)_2\right) & \text{otherwise.} \end{cases}$$

Conditioned on the event $U < p_{K,0}$, the output $\pi((K,U))_3$ is a uniform random variable; similarly for the event $U \geq p_{K,0}$. Under both outcomes, $\pi((K,U))_3$ is independent of the Poisson point process of intensity $\beta$ on $K$ given by $\pi((K,U))_2$.

As in Theorem 5.1, we input randomness from fitted shells to $\pi$. Let $\xi : F^*(X) \to [0,1]$ be defined as in Corollary 4.4. Set $\pi_{i,2} := \pi((\operatorname{cl} \mathcal{S}_i, \xi(F_i)))_2$ for $F_i \in F^*(X)$. As before, each $\pi_{i,2}$ is $P_\alpha$ almost surely a Poisson process of intensity $\beta$ on $\mathcal{S}_i$, independent of $X_{F^*(X)^c}$. Define

$$\phi(X)_2 := \sum_{F_i \in F^*(X)} \ell_i \pi_{i,2}.$$

Then $\phi(X)_2$ is a Poisson point process of intensity $\beta$ on $G$.

Set $\pi_{i,3} := \pi((\operatorname{cl} \mathcal{S}_i, \xi(F_i)))_3$. Unlike in Theorem 5.1, we now have a sequence of independent uniform random variables $\{\pi_{i,3}\}_{i \in \mathbb{N}}$, independent of $\{\pi_{i,2}\}_{i \in \mathbb{N}}$ and $X_{F^*(X)^c}$. Fix a Borel isomorphism $R : [0,1] \to A(e, 6, \rho)$, where $[0,1]$ and $A(e, 6, \rho)$ carry uniform measures. Each $\ell_i R(\pi_{i,3})$ is a Haar-uniform random point in $F_i$:

$$\ell_i R(\pi_{i,3}) \stackrel{d}{=} X|_{F_i}$$

by Remark 2.2 and Corollary 4.4. Define

$$\phi(X)_1 := X|_{F^*(X)^c} + \sum_{F_i \in F^*(X)} \ell_i R(\pi_{i,3}),$$

so by the same statements, $\phi(X)_1$ is a Poisson point process of intensity $\alpha$ on $G$.

By construction, Corollary 4.4, and Lemma 4.5, the processes $\phi(X)_1$ and $\phi(X)_2$ are independent. Set $\phi(X) := (\phi(X)_1, \phi(X)_2)$. The map $\phi$ is bijective. It is straightforward to check $\phi$ is $G$-equivariant in each coordinate, since the sets $\ell^*(X)$ and $F^*(X)$ are $G$-equivariant. Thus, $\phi$ is a factor from $\mathbb{X}_\alpha$ to $(\mathbb{X}_\alpha, \mathbb{X}_\beta)$. Each coordinate mapping is finitary by Remark 3.3 and Corollary 3.8, so $\phi$ is finitary.

The inverse map $\phi^{-1}$ serves as a finitary factor from $(\mathbb{X}_\alpha, \mathbb{X}_\beta)$ to $\mathbb{X}_\alpha$, so $\phi$ is a finitary isomorphism from $\mathbb{X}_\alpha$ to $(\mathbb{X}_\alpha, \mathbb{X}_\beta)$. □

A version of the main theorem for compactly generated groups follows directly from Proposition 5.2.

THEOREM 5.3 (Finitary isomorphism for compactly generated groups). *Let $G$ be a nondiscrete, noncompact, locally compact and compactly generated Polish group. All Poisson systems on $G$ are finitarily isomorphic.*

PROOF. Let $\alpha, \beta > 0$ and let $\mathbb{X}_\alpha$ and $\mathbb{X}_\beta$ be Poisson systems. By Proposition 5.2, each Poisson system is finitarily isomorphic to the product system $(\mathbb{X}_\alpha, \mathbb{X}_\beta)$. □

5.2. *Noncompact, compactly generated open subgroup.* Now we fix $G$ as a nondiscrete, noncompact, locally compact Polish group that is not compactly generated. Let $K$ be some compact, symmetric (so $K = K^{-1}$) neighborhood of $e \in G$. The *subgroup generated by $K$* is

$$\langle K \rangle := \bigcup_{i \geq 1} K^i,$$

and is compactly generated and $\sigma$-compact. It is simple to show $\langle K \rangle$ is open and closed (see, e.g., Proposition 2.4 in [7]) as well as of countable, possibly finite, index in $G$ ([5], Corollary 2.C.6). Further, $\langle K \rangle$ is locally compact and second countable, so $\langle K \rangle$ is uncountable Polish. It is either noncompact or compact (both instances occur; see Section 2.1 for examples). In the noncompact case, the method of coinduction applies to extend the isomorphism for compactly generated groups to $G$. We begin by realizing Stepin's method for this setting.

Fix $H$ as some noncompact, locally compact and compactly generated open subgroup of $G$. Then $H$ has countable index in $G$. We denote the set of left cosets of $H$ in $G$ as $G/H$. Since $H$ is open, it is Polish; thus, $\mathbb{M}(H)$ (the space of Borel simple point measures on $H$) and $\mathbb{M}(H)^{G/H}$ with the product topology are as well. Choose coset representatives $\{g_i\}_{i \in N} \subset G$ for $G/H$, where $N \subseteq \mathbb{N}$ enumerates the index of $H$ in $G$. Let $(\nu(g_i H))_{i \in N} \in \mathbb{M}(H)^{G/H}$.

Toward obtaining a $G$-action on $\mathbb{M}(H)^{G/H}$, define $s : G/H \to G$ so that $s$ sends a coset to its representative, and $c : G \times G/H \to H$ so that

$$c(f, gH) := s(fgH)^{-1} f s(gH)$$

for $f \in G$ and $gH \in G/H$. For $f_1, f_2 \in G$ and $gH \in G/H$, it follows from the above definition that $c(f_1 f_2, gH) = c(f_1, f_2 gH) c(f_2, gH)$ (so $c$ is a cocycle). The *action of $H \curvearrowright \mathbb{M}(H)$ coinduced to $G$* is given by

$$f(\nu(g_i H))_{i \in N} := \big(c(f^{-1}, g_i H)^{-1} \nu(f^{-1} g_i H)\big)_{i \in N}.$$

The cocycle condition can be used to show

$$f_1 f_2 (\nu(g_i H))_{i \in N} = f_1 \big(f_2 (\nu(g_i H))_{i \in N}\big),$$

so the action is well defined. Each $c(f^{-1}, g_i H)^{-1} \in H$, and recall we have a natural action $H \curvearrowright \mathbb{M}(H)$. Applying $f \in G$ to $(\nu(g_i H))_{i \in N} \in \mathbb{M}(H)^{G/H}$ permutes the sequence and shifts individual elements in the sequence by some element in $H$.

If $(\mathbb{M}(H), P_\alpha, H)$ is a Poisson system over $H$, then

$$\mathbb{X}_\alpha^H := (\mathbb{M}(H)^{G/H}, P_\alpha^{G/H}, G)$$

is the *Poisson system coinduced from* $H \curvearrowright (\mathbb{M}(H), P_\alpha)$. With coinduction, we prove some systems of interest are isomorphic in Lemma 5.4. Factor maps between $(\mathbb{M}(G), P_\alpha, G)$ and $\mathbb{X}_\alpha^H$ are defined as expected. Such factors are finitary if each of their coordinate mappings are.

LEMMA 5.4 (Coset copy isomorphism). *Let $(\mathbb{M}(H), P_\alpha, H)$ be a Poisson system over $H$. The Poisson system coinduced from $H \curvearrowright (\mathbb{M}(H), P_\alpha)$ is finitarily isomorphic to $(\mathbb{M}(G), P_\alpha, G)$.*

PROOF. Let $(\nu(g_i H))_{i \in N} \in \mathbb{M}(H)^{G/H}$. Define $\psi : \mathbb{M}(H)^{G/H} \to \mathbb{M}(G)$ by

$$\psi\big((\nu(g_i H))_{i \in N}\big) := \sum_{i \in N} g_i \nu(g_i H).$$

The map $\psi$ takes a sequence in $\mathbb{M}(H)^{G/H}$ and sends each element of the sequence to the copy of $H$ in $G$ identified by the element's index. If each element in the input sequence were a Poisson point process on $H$ with intensity $\alpha$, then the output would be a Poisson point process on $G$ with intensity $\alpha$ by Remark 3.5, since $G/H$ is countable. The map $\psi$ is measurable, bijective and finitary with finitary inverse.

To show $\psi$ is a finitary isomorphism, it remains to check $\psi$ and its inverse are $G$-equivariant on a set of full measure. In fact, $\psi$ is $G$-equivariant on $\mathbb{M}(H)^{G/H}$. Let $f \in G$. We have

$$f\psi\big((\nu(g_i H))_{i \in N}\big) = f \sum_{i \in N} g_i \nu(g_i H)$$

as well as

$$\psi\big(f(\nu(g_i H))_{i \in N}\big) = \psi\big((c(f^{-1}, g_i H)^{-1} \nu(f^{-1} g_i H))_{i \in N}\big),$$

where for each $i$, $(c(f^{-1}, g_i H)^{-1} = h_i^{-1}$ for some $h_i^{-1} \in H$ and $g_{k_i} \in \{g_i\}$ such that $f^{-1} g_i = g_{k_i} h_i$. So, we may write

$$\psi\big(f(\nu(g_i H))_{i \in N}\big) = \psi\big((\nu(f^{-1} g_i H) \circ h_i)_{i \in N}\big)$$

$$= \sum_{i \in N} g_i \big(\nu(f^{-1} g_i H) \circ h_i\big)$$

$$= \sum_{i \in N} f g_{k_i} h_i \big(\nu(g_{k_i} H) \circ h_i\big)$$

$$= f \sum_{i \in N} g_{k_i} \nu(g_{k_i} H) = f\psi\big((\nu(g_i H))_{i \in N}\big).$$

A similar argument shows $\psi^{-1}$ is $G$-equivariant. Hence, $\mathbb{X}_\alpha^H$ is finitarily isomorphic to $(\mathbb{M}(G), P_\alpha, G)$. □

Suppose we have a factor map $\phi : (\mathbb{M}(H), P_\alpha, H) \to (\mathbb{M}(H), P_\beta, H)$, for some $\alpha, \beta > 0$. Keeping the same notation for $(\nu(g_i H))_{i \in N} \in \mathbb{M}(H)^{G/H}$ as before, define $\Psi : \mathbb{M}(H)^{G/H} \to \mathbb{M}(H)^{G/H}$ by

$$\Psi\big((\nu(g_i H))_{i \in \mathbb{N}}\big) := \big(\phi(\nu(g_i H))\big)_{i \in \mathbb{N}}.$$

We say $\Psi$ is the *factor map coinduced by* $\phi$. Such a factor is finitary if each coordinate mapping is.

LEMMA 5.5 (Coinduced isomorphism). *Let $\alpha, \beta > 0$ and let $(\mathbb{M}(H), P_\alpha, H)$ and $(\mathbb{M}(H), P_\beta, H)$ be Poisson systems over $H$. The Poisson systems coinduced from $H \curvearrowright (\mathbb{M}(H), P_\alpha)$ and $H \curvearrowright (\mathbb{M}(H), P_\beta)$ are finitarily isomorphic.*

PROOF. By Theorem 5.3, there exists a finitary isomorphism from $(\mathbb{M}(H), P_\alpha, H)$ to $(\mathbb{M}(H), P_\beta, H)$. Let $\phi$ be such an isomorphism. Let $\Psi$ be the factor map coinduced by $\phi$. Then $\Psi$ is measurable, bijective, finitary with finitary inverse and by Remark 3.5,
$$P_\alpha^{G/H} \circ \Psi^{-1} = P_\beta^{G/H}.$$

Let $f \in G$ and $(\nu(g_i H))_{i \in N}$. With $f$ acting as in Lemma 5.4, here we have
$$f \Psi(\nu(g_i H)) = f \phi(\nu(g_i H)) = \phi(\nu(f^{-1} g_i H) \circ h_i)$$
for some $h_i \in H$, since $\phi$ is $H$-equivariant and $\phi(\nu(g_i H))$ is an element of $\mathbb{M}(H)$ indexed by $G/H$, so the $f$-action is well defined. On the other hand,
$$\Psi(f \nu(g_i H)) = \Psi(\nu(f^{-1} g_i H) \circ h_i) = \phi(\nu(f^{-1} g_i H) \circ h_i).$$
Thus, $\Psi$ is a finitary isomorphism between the Poisson systems coinduced from $H \curvearrowright (\mathbb{M}(H), P_\alpha)$ and $H \curvearrowright (\mathbb{M}(H), P_\beta)$. □

Theorem 5.6 follows from Lemmas 5.4 and 5.5.

THEOREM 5.6 (Finitary isomorphism for noncompact subgroup case). *Let $G$ be a nondiscrete, noncompact, locally compact Polish group with a noncompact, compactly generated open subgroup $H$. All Poisson systems on $G$ are finitarily isomorphic.*

PROOF. By Lemma 5.5, the Poisson systems coinduced from $H \curvearrowright (\mathbb{M}(H), P_\alpha)$ and $H \curvearrowright (\mathbb{M}(H), P_\beta)$ are finitarily isomorphic, and by Lemma 5.4, they are finitarily isomorphic to $(\mathbb{M}(G), P_\alpha, G)$ and $(\mathbb{M}(G), P_\beta, G)$, respectively. Thus, the systems $(\mathbb{M}(G), P_\alpha, G)$ and $(\mathbb{M}(G), P_\beta, G)$ are finitarily isomorphic. □

Although we do not need a version of Proposition 1.5 to prove Theorem 5.6, we do need one to eventually prove Corollary 1.6. We state this version in Proposition 5.7, which follows from Theorem 5.2 and Lemma 5.4.

PROPOSITION 5.7. *Let $G$ be a nondiscrete, noncompact, locally compact Polish group with a noncompact, compactly generated open subgroup $H$. Let $\alpha, \beta > 0$. There exists a finitary isomorphism from $\mathbb{X}_\alpha$ to $(\mathbb{X}_\alpha, \mathbb{X}_\beta)$.*

**6. General case.** In Section 6.1, we assume $G$ is nondiscrete, noncompact, locally compact Polish and does not contain any noncompact, compactly generated open subgroups. It is easy to see by the argument at the start of Section 5.2 that $G$ must contain a compact open subgroup. (In fact, $G$ contains infinitely many compact open subgroups; see Remark 6.3.) Let $H$ be a compact open subgroup of $G$, although we will define $H$ a bit more specifically shortly. Our goal is to define markers on $H$ and cosets of $H$ and a Voronoi tessellation on $G$ via its subgroup structure.

Once we have proved Theorem 1.1 for this case, we will prove the general case in Section 6.2.

6.1. *Compact subgroup chain.* Unlike in Section 5.2, we cannot utilize prior results for noncompact, compactly generated groups. Instead, we modify the results in Sections 3 and 4 to work for our current needs. We again appeal to a word metric to define markers on $H$. The subgroup $H$ is trivially compactly generated, but assigning $H$ as its own generating set yields a rather useless word metric. We prove Lemma 6.1 to find a generating set so that $H$ is spacious enough under the corresponding word metric to allow marker independence properties.

We do not know if Lemma 6.1 can be adjusted to allow such a generating set, on the other hand, to be spacious enough with respect to Haar measure so that fitted shells may host the Poisson coupling from Lemma 3.1. For this reason, we do not construct a splitting map for Poisson systems over $G$ as defined in this section (see Section 7).

LEMMA 6.1. *Let $G$ be a noncompact, locally compact Polish group without a noncompact, compactly generated open subgroup and fix $n \in \mathbb{N} \setminus \{0\}$. Then $G$ contains a compact, symmetric neighborhood $S$ of $e$ such that $S^n \neq S^{n+1}$.*

PROOF. Let $V$ be a compact, symmetric neighborhood of $e$ in $G$, and set $J := \langle V \rangle$. Certainly $J$ is a compactly generated open subgroup of $G$, so by the assumptions on $G$, it must be compact. Fix $m \in \mathbb{N} \setminus \{0\}$, and let $W$ be another compact, symmetric neighborhood of $e$ with $V \subset W$ such that $J$ has finite index at least $m$ in the compact open subgroup $L := \langle W \rangle$. Such a $W$ must exist since $G$ is noncompact without any noncompact, compactly generated open subgroups.

We want a subgroup of $J$ that is normal in $L$. Set

$$K := \bigcap_{g \in L} gJg^{-1} \leq J.$$

Then $K \triangleleft L$. Since $J$ has finite index in $L$, the normalizer subgroup

$$N_L(J) := \{g \in L : gJg^{-1} = J\}$$

of $J$ in $L$ has finite index in $L$. The index of $N_L(J)$ in $L$ is the number of conjugacy classes of $J$ in $L$, so the intersection in the definition of $K$ is finite, implying $K$ is open, and hence compact.

Consider $F := L/K$, a finite group. Any finite group $D$ has a minimal, symmetric generating set $R$ with order at most $\log_2 |D|$. For large enough $|D|$, the final inequality in the chain of inequalities

$$|R^n| \leq |R|^n \leq (\log_2 |D|)^n < |D|$$

holds. We have a group $F$ such that $m \leq |F|$ for any $m \in \mathbb{N} \setminus \{0\}$. Let $T$ be a minimal, symmetric generating set for $F$ with $|T| \leq \log_2 |F|$. Thus, we can find $F$ and $T$ such that $|T^n| < |F|$, meaning $T^n \subsetneq T^{n+1}$.

If $T$ generates $F$ so that $T^n \subsetneq T^{n+1}$, then $T \cup K$ generates $L$ and $(T \cup K)^n \subsetneq (T \cup K)^{n+1}$. Set $S := T \cup K$. Note $S$ is a compact, symmetric neighborhood of $e$. □

Now fix $S$ as a compact, symmetric neighborhood of $e$ in $G$ such that $S^{60} \neq S^{61}$, and set $H := \langle S \rangle$, so $H$ is a compact open subgroup of $G$. We proceed to redefine marker tools from Section 3. Although we may use the same terminology, some definitions are different. Since we do not use the markers from this section to construct a splitting map, we define shell values explicitly rather than in terms of $\rho$ (if we defined them in terms of $\rho$, they would not be guaranteed to fit inside of $H$).

Let $d_S : H \times H \to \mathbb{N}$ be the word metric defined by $S$ on $H$, so that for any $g, h \in H$,

$$d_S(g, h) = \min \left\{ n \geq 0 : \begin{array}{l} \text{there exist } s_1, \ldots, s_n \in S \\ \text{such that } g^{-1}h = s_1 \ldots s_n \end{array} \right\}.$$

Note the diameter of $H$ with respect to $d_S$, denoted as $\operatorname{diam} H$, is finite and at least 60. Balls are again defined to contain their boundary. Let $X$ be a Poisson point process on $H$ with intensity $\alpha$. We say $h \in H$ is a *seed* if $B(h, \operatorname{diam} H) = H$ satisfies the following under $X$:

1. For every ball $B$ of radius 1 such that $B \subset A(h, 16, 20)$, $N(B) \geq 1$.
2. The shell $A(h, 20, \operatorname{diam} H)$ is empty.

With nonzero probability, $H$ does not contain any seeds. However, we are ultimately interested in the collection of cosets of $H$ covering $G$. Because $H$ is compact, the index of $H$ in $G$ is countably infinite.

The proof of Lemma 6.2 follows the same argument as the proof of Lemma 3.2. However, in Lemma 6.2, we gain the additional property that if seeds exist in $H$ under $X$, they must belong to the same cluster (to be defined as an equivalence class).

LEMMA 6.2 (Seed distances for a compact group). *Let $X$ be a Poisson point process on $H$ with intensity $\alpha$, and let $g, h \in H$ be seeds under $X$. Then $d_S(g, h) < 2$.*

We define an equivalence relation by setting seeds $g$ and $h$ equivalent whenever $d_S(g, h) < 2$. Given a seed $g$, its equivalence class of seeds is $[g]$. The *core* of $[g]$ is

$$\mathcal{C}[g] := \bigcap_{h \in [g]} B(h, 3).$$

Cores containing a unique point of $X$ are *identifiable*, and the unique $X$-point in an identifiable core is its *landmark* $\ell[g] := X|_{\mathcal{C}[g]}$. The *fitted shell* of a landmark is

$$F[g] := A(\ell[g], 6, 10).$$

Lemma 6.2 implies that Lemma 3.4 applies here, with $\mu \in \mathbb{M}(H)$ such that $\mu$ contains a landmark (necessarily unique in $H$). If $X|_{F[g]} = 1$, the associated seed $g$ with its class, core, landmark, and fitted shell are *harvestable*, and $X|_{F[g]}$ is their *harvest*.

Let $\{g_j\}_{j \in \mathbb{N}}$ be coset representatives of $G/H$ with fixed enumeration, so

$$G = \bigsqcup_{j \in \mathbb{N}} g_j H.$$

Definitions for seeds, cores, et cetera hold not just on $H$, but on any coset of $H$ as well: If $g$ and $h$ belong to the same coset of $H$, then $g^{-1}h \in H$ and $d_S(g, h) = d_S(e, g^{-1}h)$ is well defined. A coset may contain at most one landmark, by the definition of a seed and Lemma 6.2.

Let $X$ be a Poisson point process on $G$ with intensity $\alpha$. Borel–Cantelli implies infinitely many cosets of $H$ contain a landmark, since the index of $H$ in $G$ is infinite. Denote the set of harvestable landmarks under $X$ as $\ell^*(X) := \{\ell_i\}_{i \in \mathbb{N}}$. We use $\ell^*(X)$ to define a Voronoi tessellation on $G$. We cannot use the definition for restricted Voronoi cells, however, and neither can we modify our coinduction argument. We no longer have a word metric on $G$ defined by a compact set, and in the case for $H$ noncompact, restricted Voronoi cells lived in $H$, but cells must live in $G$ here. If we were to replace the word metric on $G$ with some proper, left-invariant, compatible metric, we would have no control over the Haar measure of the topological boundary of a cell. We take a new approach, which relies on a structural property of $G$, stated below; see [5], Proposition 2.C.3 for a proof.

REMARK 6.3 (Compact subgroup chain). Let $G$ be a locally compact and $\sigma$-compact group. Then there exists a monotonically increasing sequence $\{H_n\}_{n\geq 1}$ of compactly generated open subgroups of $G$ such that

$$G = \bigcup_{n\geq 1} H_n.$$

If $G$ does not contain any compactly generated, noncompact open subgroups, then additionally each $H_n$ is compact, and the sequence is infinite.

We define a metric on $G/H$ to replace the word metric on $G$ from Section 3. Let $\{H_n\}_{n\geq 1}$ be a sequence of compact open subgroups of $G$ as described in Remark 6.3, with $H_0 := H \subseteq H_1$, and let $gH, fH \in G/H$. Define $d_G : G/H \times G/H \to \mathbb{N}$ so that $d_G(gH, fH) := n$, where $n$ is the minimum value such that $g^{-1}f \in H_n$. It is straightforward to check $d_G$ is a metric on $G/H$, and is in fact an ultrametric. We use $d_G$ to define a Voronoi tessellation on $G$. First, we need a tie-breaking map.

Enumerate harvestable cores under $X$ as $\{\mathcal{C}_i\}_{i\in\mathbb{N}}$. Let $\mathcal{A}$, $d'$ and the equivalence relation $\sim$ on $\mathcal{A}$ be defined as in Section 3, except we use $d_S$ as defined in this section. Fix a Borel isomorphism $T : \mathcal{A}/\sim \to [0, 1]$, and set

$$T_i := T(\ell_i^{-1}\mathcal{C}_i, e)$$

for all $\ell_i \in \ell^*(X)$. Since the new Voronoi cells are quite different from those in Section 3, we give them a new name.

The *corestricted Voronoi cell* of $\ell_i \in \ell^*(X)$ is the set

$$\mathcal{V}^*(\ell_i) := \left\{ \bigcup fH \subseteq G : \begin{array}{l} \text{for each } \ell_k \in \ell^*(X), \text{ either } d_G(\ell_i H, fH) < d_G(\ell_k H, fH) \text{ or} \\ d_G(\ell_i H, fH) = d_G(\ell_k H, fH) \text{ and } T_i \leq T_k \end{array} \right\}$$

and the *corestricted Voronoi tessellation* of $X$ is

$$\mathcal{V}^*(X, H) := \{\mathcal{V}^*(\ell_i)\}_{i\in\mathbb{N}}.$$

Corestricted Voronoi tessellations are measurable partitions of $G$ almost surely, and are $G$-equivariant since $d_G$ is left-invariant. We note that although $d_G$ suffices to define a Voronoi tessellation on $G$, we cannot use it to define markers. Ultrametrics are not conducive to proving the property in Lemma 6.2, while word metrics are.

In Lemma 6.4 and Corollary 6.5, we prove corestricted Voronoi cells have similar properties to restricted Voronoi cells, in analogues to Lemma 3.7, Corollary 3.8 and Lemma 3.10. For $g \in G$, we denote the unique harvestable landmark in the same corestricted Voronoi cell as $g$ under $X$ as $\ell(X, g)$, and the corestricted Voronoi cell itself as $\mathcal{V}(X, g)$.

LEMMA 6.4 (Corestricted Voronoi cells have finite expected measure). *Let $X$ be a Poisson point process on $G$ with intensity $\alpha$ and corestricted Voronoi tessellation $\mathcal{V}^*(X, H)$. Then*

$$\mathbb{E}\big(\lambda(\mathcal{V}(X, e))\mathbb{1}_{\ell(X,e)\in H}\big)$$

*is finite $P_\alpha$ almost surely.*

PROOF. Our approach is almost the same as in Lemma 3.7. We reuse the map $M$ and slightly modify the map $\Lambda$. Let $M : \mathbb{M} \times G \to \mathbb{M} \times \mathbb{M}$ be defined as in Lemma 3.7. Note that we have $\ell(\mu, e) \in \ell(\mu, e)H$. Let $\mu, \nu \in \mathbb{M}$ and $g \in G$. Define $\Lambda : \mathbb{M} \times \mathbb{M} \to \{0, 1\}$ so that

$$\Lambda(\mu, \nu) := \begin{cases} 1 & \text{if } \nu = g\mu \text{ for some } g \text{ such that } g^{-1} \in \ell(\mu, e)H, \\ 0 & \text{otherwise.} \end{cases}$$

Then
$$\int_G \Lambda(M(\mu, g)) \, d\lambda(g) = \lambda(\ell(\mu, e)H) = \lambda(H).$$

As before, integrating the above over $\mathbb{M}$ with respect to the probability measure $P_\alpha$ yields the same result, and since $\lambda$ is $\sigma$-finite and left-invariant, we have

$$\int_\mathbb{M} \int_G \Lambda(M(\mu, g)) \, d\lambda(g) \, dP_\alpha(\mu) = \int_\mathbb{M} \int_G \Lambda(M(g^{-1}\mu, g)) \, d\lambda(g) \, dP_\alpha(\mu),$$

where the last integral is
$$E := \mathbb{E}(\lambda(\mathcal{V}(\mu, e))\mathbb{1}_{\ell(\mu,e)\in H}),$$

because $\Lambda(M(g^{-1}\mu, g)) = 1$ if and only if $\ell(\mu, g) \in H$ and $g \in \mathcal{V}(\mu, e)$. Therefore, $E = \lambda(H) < \infty$. □

COROLLARY 6.5 (Compact corestricted cells). *Let $X$ be a Poisson point process on $G$ with intensity $\alpha$ and corestricted Voronoi tessellation $\mathcal{V}^*(X, H)$. Then each corestricted Voronoi cell is $P_\alpha$ almost surely compact.*

PROOF. By definition, a corestricted Voronoi cell is a union of cosets of $H$. Suppose for a contradiction that, with positive probability under $P_\alpha$, there exists a corestricted Voronoi cell $V \in \mathcal{V}^*(X, H)$ such that

$$V = \bigsqcup_{j \in J} g_j H,$$

where $J$ is some infinite subset of $\mathbb{N}$. But then $\lambda(V) = \infty$, so $\mathbb{E}(\lambda(\mathcal{V}(\mu, e))\mathbb{1}_{\ell(\mu,e)\in H}) = \infty$. This contradicts Lemma 6.4.

We have shown each corestricted Voronoi cell is $P_\alpha$ almost surely a finite union of cosets of $H$, so each cell is $P_\alpha$ almost surely compact. □

As before, harvestable fitted shells are sources of randomness for populating output processes on Voronoi cells. Enumerate the set of harvestable fitted shells under $X$ as $\{F_i\}_{i \in \mathbb{N}}$, and let $F^*(X)$ be the union of harvestable fitted shells under $X$. Lemma 4.2 holds on each $H$ copy with $G = H = A$ by Lemma 6.2, and cosets of $H$ are disjoint, so we have independence between coset copies. Thus, Proposition 4.1 and Corollaries 4.3 and 4.4 hold for $G$ and $F^*(X)$ as defined in this section. We redefine the map $\xi$ in Corollary 4.4 for this section; it remains the same as before, except the map $R$ in its definition has domain $A(e, 6, 10)$.

THEOREM 6.6 (Finitary factor for remaining groups). *Let $G$ be a nondiscrete, noncompact, locally compact Polish group without any noncompact, compactly generated open subgroups. Let $\alpha, \beta > 0$. There exists a finitary factor from $\mathbb{X}_\alpha$ to $\mathbb{X}_\beta$.*

The proof of Theorem 6.6 follows the precise argument of Theorem 5.1. Corestricted Voronoi tessellations replace restricted Voronoi tessellations. Corestricted Voronoi cells are almost surely compact by Corollary 6.5, so the map $\pi$ applies, and the constructed factor is finitary.

The corestricted Voronoi tessellation of $X$ depends on $X|_{F^*(X)^c}$, the fixed enumeration $\{g_i\}_{i \in \mathbb{N}}$ of coset representatives, and the fixed tie-breaking map $T$, so it is independent of $X|_{F^*(X)}$. Thus, Lemma 4.5 holds in this setting, so Proposition 5.2 holds here as well, which implies Proposition 1.5.

6.2. *The factor map and isomorphisms.* Finally, fix $G$ as a nondiscrete, noncompact, locally compact Polish group. Theorem 1.2 follows from Theorems 5.1, 5.6 and 6.6. Proposition 1.5 follows from Propositions 5.2 and 5.7 and Theorem 6.6. Theorem 1.1 follows from Theorems 5.3, 5.6 and Proposition 1.5. The constructions used to prove Theorems 1.1 and 1.2 imply Theorem 1.3.

PROOF OF COROLLARY 1.6. Let $\alpha, \beta, \gamma > 0$ and let $X$, $Y$ and $Z$ be Poisson point processes on $G$ with intensities $\alpha$, $\beta$ and $\gamma$, respectively. By Proposition 1.5, the Poisson system $\mathbb{X}_\alpha$ is isomorphic to $(\mathbb{X}_\alpha, \mathbb{X}_\beta)$ and $\mathbb{X}_\beta$ is isomorphic to $(\mathbb{X}_\beta, \mathbb{X}_\gamma)$. Theorem 1.1 implies $\mathbb{X}_\alpha$ is isomorphic to $\mathbb{X}_\beta$, and thus $(\mathbb{X}_\beta, \mathbb{X}_\gamma)$. □

**7. Poisson splitting.** We conclude the paper with a proof of Theorem 1.4 and two questions regarding extensions of the theorem. Once again, we let $G$ be a nondiscrete, noncompact, locally compact, compactly generated Polish group. We first prove the theorem for this case, and then extend via coinduction to the case when the group is not compactly generated, but contains a noncompact, compactly generated open subgroup.

Let $\alpha > \beta > 0$. We construct a monotone, finitary factor $\Phi$ from $\mathbb{X}_\alpha$ to $\mathbb{X}_\beta$ with the additional property that if $X$ is a Poisson point process on $G$ with intensity $\alpha$, then $X - \Phi(X)$ is a Poisson point process on $G$ with intensity $\alpha - \beta$. We follow many of the same steps as in [8]. Randomness to power the splitting map comes from fitted shells under the input process $X$, which contain precisely either one or two points of $X$. Within these fitted shells, we deterministically retain points from those with one point (1-fitted shells) and delete points from those with two points (2-fitted shells). This determinism frees up the randomness within these shells to apply to points of $X$ elsewhere. In particular, we apply a splitting based on the coupling in Lemma 3.1 to all fitted shells, and a splitting based on Remark 7.2 to ambient points of $X$ outside of fitted shells.

For the splitting map to delete $X$-points $G$-equivariantly, we place two Voronoi tessellations on $G$, one with landmarks of 1-fitted shells as sites, and the other with landmarks of 2-fitted shells as sites. Let $d_S$ be the word metric on $G$ defined in Section 3. Note 1-fitted shells are harvestable fitted shells, so the restricted Voronoi tessellation $\mathcal{V}^*(X)$ is the first tessellation described. We also refer to $\mathcal{V}^*(X)$ as the 1-Voronoi tessellation of $X$. We denote the set of landmarks associated to 2-fitted shells as $\ell^{**}(X)$ and the union of 2-fitted shells as $F^{**}(X)$, and we set $\ell'(X) := \ell^*(X) \cup \ell^{**}(X)$ and $F'(X) := F^*(X) \cup F^{**}(X)$.

Let $T$ be the tie-breaking map from Section 3. The 2-*Voronoi cell* of the 2-landmark $\ell[g] \in \ell^{**}(X)$ is

$$\mathcal{V}^{**}(\ell[g]) := \left\{ f \in G : \begin{array}{l} \text{for each } \ell[h] \in \ell^{**}(X), \text{ either } d_S(\ell[g], f) < d_S(\ell[h], f) \text{ or} \\ d_S(\ell[g], f) = d_S(\ell[h], f) \text{ and } T_{C[g], \ell[g]} \leq T_{C[h], \ell[h]} \end{array} \right\},$$

and the 2-*Voronoi tessellation* of $X$ is $\mathcal{V}^{**}(X) := \{\mathcal{V}^{**}(\ell[g])\}_{\ell[g] \in \ell^{**}(X)}$. We have that 2-Voronoi cells are almost surely precompact with measure zero boundary by the same argument as for 1-Voronoi cells (Corollary 3.8 and Lemma 3.10). Fix a Borel isomorphism $Q : G \to [0, 1]$. We use $Q$ to enumerate points in 1- and 2-Voronoi cells, after shifting cells to $e \in G$ via their sites.

Over $\mathbb{R}^n$, distances between points of a Poisson point process are almost surely unique. Holroyd, Lyons and Soo exploit this in [8] to equivariantly link randomness inside their versions of 1- and 2-fitted shells to points outside such shells. We use 1- and 2-Voronoi tessellations and $Q$ instead. Any $g \in G$ will almost surely belong to unique 1- and 2-Voronoi cells. If $g$ is "turned on" under $X$, and $g \notin F'(X)$, then we decide whether to delete $g$ using information from its associated 1- and 2-fitted shells. We queue $X$-points to receive information from their relevant fitted shells via the $Q$ enumeration.

Proposition 4.1 holds for $F^{**}(X)$ by the same argument as for the original $F^*(X)$. Corollary 7.1 is based on Corollary 4.4 but contains slight changes. We enumerate sets of 2-objects as $\{C_i^{**}\}_{i\in\mathbb{N}}$, $\{\ell_i^{**}\}_{i\in\mathbb{N}}$, $\{F_i^{**}\}_{i\in\mathbb{N}}$ and $\{\mathcal{V}_i^{**}\}_{i\in\mathbb{N}}$ for harvestable 2-cores, 2-landmarks, 2-fitted shells and 2-Voronoi cells, respectively.

COROLLARY 7.1 (Source independence for splitting). *Let $X$ be a Poisson point process on $G$ with intensity $\alpha$, and conditioned on $F'(X)$, let*

$$U := \{U_{F_i'}\}_{F_i' \in F'(X)}$$

*be a sequence of independent uniform random variables, independent of $X$. There exists a measurable map $\xi : F'(X) \to [0,1]$ such that*

$$(X|_{F'(X)^c}, F'(X), \{\xi(F_i')\}_{F_i' \in F'(X)}) \stackrel{d}{=} (X|_{F'(X)^c}, F'(X), U).$$

PROOF. By Proposition 4.1, we have that

$$\mathcal{S}'(X) := \{(\ell_i)^{-1} X|_{F_i}\}_{i\in\mathbb{N}} \cup \{(\ell_i^{**})^{-1} X|_{F_i^{**}}\}_{i\in\mathbb{N}}$$

is a sequence of independent Poisson point processes on $A := A(e, 6, \rho)$ independent of $(X|_{F'(X)^c}, F'(X))$. Let $V$ be a uniform random variable on $A$, independent of $X$. Fix a Borel isomorphism $R : A \to [0,1]$ so that $R(V)$ is a uniform random variable on $[0,1]$. Denote the shifted process point of $(\ell_i)^{-1} X|_{F_i}$ as $v_i$ and the shifted process points of $(\ell_i^{**})^{-1} X|_{F_i^{**}}$ as $v_{i,1}$ and $v_{i,2}$, for $\ell_i \in \ell^*(X)$ and $\ell_i^{**} \in \ell^{**}(X)$.

We remark that if $U_1$ and $U_2$ are independent uniform random variables, then $U_1 \oplus U_2$ is a uniform random variable independent of $U_1$ and $U_2$, where $\oplus$ denotes addition modulo 1 (see [8], Lemma 9). Let $F_i' \in F'(X)$. Define $\xi : F'(X) \to [0,1]$ by

$$\xi(F_i') := \begin{cases} R(v_{i,1}) & \text{if } F_i' \in F^*(X), \\ R(v_{i,1}) \oplus R(v_{i,2}) & \text{otherwise.} \end{cases}$$

Each $\xi(F_i')$ for $F_i' \in F'(X)$ is uniformly distributed on $[0,1]$, and the sequence of independent uniform random variables $\{\xi(F_i')\}_{F_i' \in F'(X)}$ is independent of $(X|_{F'(X)^c}, F'(X))$, which implies the statement. $\square$

We want to construct a factor map $\Phi$ from $\mathbb{X}_\alpha$ to $\mathbb{X}_\beta$ with the additional property that, if $X$ is a Poisson point process with intensity $\alpha$, then $X - \Phi(X)$ is a Poisson process with intensity $\alpha - \beta$. This map should delete $X$-points with randomness from $X$ itself in such a way as to end up with a Poisson point process of intensity $\beta$.

Given a single uniform random variable, we can apply the following map to any $X$-point. Let $(g, x) \in G \times [0,1]$. Define $D : G \times [0,1] \to \mathbb{M}$ by $D(g, x) = \mathbb{1}_{x \leq \beta/\alpha} \delta(g)$. We may extend this map to all points of $X$ (but not $G$-equivariantly). Recall our Borel isomorphism $Q : G \to [0,1]$. Let $\mu \in \mathbb{M}$, and enumerate the process points of $\mu$ as $\{x_i\}_{i\in\mathbb{N}}$, where $Q(x_i) < Q(x_{i+1})$ for almost every $i \in \mathbb{N}$. Recall the bit map $b$ from Proposition 5.2. We modify $b$ to split $x \in [0,1]$ into infinitely countably many points in $[0,1]$; for the remainder of the paper, we consider $b : [0,1] \to [0,1]^{\mathbb{N}\setminus\{0\}}$. Denote the output of $b(x)$ as $\{b(x)_i\}_{i\geq 1}$. Now, define $D' : \mathbb{M} \times [0,1] \to \mathbb{M}$ so that

$$D'(\mu, x) = \sum_{i\in\mathbb{N}} \mathbb{1}_{b(x)_{i+1} \leq \beta/\alpha} \delta(x_i).$$

The map $D'$ deletes points of $\mu$ independently with probability $\beta/\alpha$, so inputting a Poisson point process on $G$ with intensity $\alpha$ to $D'$ along with a uniform random variable (and a Borel isomorphism $Q$) yields a Poisson point process with intensity $\beta$.

Certainly $D'$ is not bijective. Since it is impossible for a deletion map to be bijective, we do not concern ourselves with the loss of information here.

REMARK 7.2 (Independent splitting). Suppose $X$ is a Poisson point process on $G$ with intensity $\alpha$ and $U$ is a uniform random variable independent of $X$. Let $A \in \mathcal{B}(G)$. Then $D'(X|_A, U)$ and $X|_A - D'(X|_A, U)$ are independent Poisson point processes on $A$ with intensities $\beta$ and $\alpha - \beta$, respectively.

Remark 7.2 appears as Lemma 19 in [8] and is immediate from the definition of $D'$. To obtain a splitting factor map (the splitting in Remark 7.2 is not a factor), we apply the independent splitting map $D$ to $X$-points outside of $F'(X)$. Inside harvestable fitted shells, we apply a map based on Lemma 3.1. This map always keeps $X$-points in 1-fitted shells and always deletes $X$-points in 2-fitted shells. We consider all harvestable fitted shells and not just 1- and 2-fitted shells so that we may transfer mass to output the desired coupling. We first construct the Lemma 3.1 map on $A(e, 6, \rho)$.

Proposition 7.3 is based on Proposition 8 in [8]. The proofs are almost the same. For any set $A$, we denote the set of all subsets of $A$ as $\mathcal{P}(A)$.

PROPOSITION 7.3 (Splitting on a fitted shell). *Let $X$ be a Poisson point process on $A := A(e, 6, \rho)$ with intensity $\alpha$, and let $U$ be a uniform random variable independent of $X$. There exists a map $S : \mathbb{M}(A) \times [0, 1] \to \mathbb{M}(A)$ with the following properties*:

1. $S(X, U)$ *and* $X - S(X, U)$ *are Poisson point processes on $A$ with intensities $\beta$ and $\alpha - \beta$, respectively.*
2. *If $g \in A$, then $S(X, U)|_{\{g\}} \leq X|_{\{g\}}$.*
3. *If $N(X) = 1$, then $S(X, U) = X$, and if $N(X) = 2$, then $S(X, U) = \varnothing$.*

PROOF. Let $\mu \in \mathbb{M}(A)$. We define $S$ piecewise with outputs depending on the number of points of $\mu$. Fix a Borel isomorphism $R : A \to [0, 1]$, so that we have an ordering on $A$ almost surely. For $k \in \mathbb{N} \setminus \{0\}$, we denote the set $\{1, \ldots, k\}$ as $[k]$. Let $j \in \mathbb{N}$ with $j \leq k$ and define $s_{j,k} : [0, 1] \to \mathcal{P}([k])$ so that $s_{j,k}(U)$ is uniformly distributed over subsets of size $j$ of $[k]$.

Define $\iota : X \to [N(X)]$ so that for $g, h \in X$, we have $\iota(g) < \iota(h)$ if and only if $R(g) < R(h)$ (here we abuse notation and think of $g, h$ both as $X$-points and elements in $A$). Note $N(X)$ is a Poisson random variable with mean $\alpha \cdot \lambda(A)$, and we have defined $A$ (in Section 3) so that $\lambda(A) > C_r/\alpha$. By Lemma 3.1, there exists a Poisson random variable $Z$ with mean $r\alpha \cdot \lambda(A) = \beta \cdot \lambda(A)$ such that

$$\mathbb{P}(N(X) - Z = 0 | N(X) = 1) = \mathbb{P}(Z = 0 | N(X) = 2) = 1$$

and $Y := N(X) - Z$ is a Poisson random variable with mean $(\alpha - \beta) \cdot \lambda(A)$.

Let $b$ be the bit map referenced in this section (or the one from Proposition 5.2), so that $b(U)_1$ and $b(U)_2$ are independent uniform random variables. We use $b(U)_1$ to choose the number of $X$-points to keep whenever $N(X) > 2$. Define $\mathcal{F}$ so that

$$(N(X), \mathcal{F}(N(X), U)) \stackrel{d}{=} (N(X), Z).$$

Now we define $S$ dependent on the value of $(\mathcal{F}(N(X), b(U)_1)$. Suppose $N(X) = k$ and $\mathcal{F}(k, b(U)_1) = j$. We use $b(U)_2$ to choose $j$ points uniformly at random to retain from the $k$ points of $X$. Set

$$S(X, U) := \sum_{g \in X : \iota(g) \in s_{j,k}(b(U)_2)} \delta(g) \quad \text{whenever } N(X) = k \text{ and } \mathcal{F}(k, b(U)_1) = j.$$

Then $S$ satisfies properties (1), (2) and (3). □

THEOREM 7.4 (Splitting factor for compactly generated groups). *Let G be a nondiscrete, noncompact, locally compact and compactly generated Polish group. For all $\alpha > \beta > 0$, there exists a monotone, finitary factor $\Phi : \mathbb{X}_\alpha \to \mathbb{X}_\beta$ with the additional property that if X is a Poisson point process on G with intensity $\alpha$, then $X - \Phi(X)$ is a Poisson point process on G with intensity $\alpha - \beta$.*

PROOF. Let $X$ be a Poisson point process on $G$ with intensity $\alpha$, and let $\xi : F'(X) \to [0, 1]$ be defined as in Corollary 7.1. Then $U := \{\xi(F'_i)\}_{F'_i \in F'(X)}$ is a sequence of independent uniform random variables independent of $(X|_{F'(X)^c}, F'(X))$. Recall

$$b : [0, 1] \to [0, 1]^{\mathbb{N} \setminus \{0\}}$$

is the bit splitting map. For each $\xi(F'_i) \in U$, the sequence $\{b(\xi(F'_i))_j\}_{j \geq 1}$ is also composed of independent uniform random variables and is independent of $(X|_{F'(X)^c}, F'(X))$.

Define $S : \mathbb{M}(A) \times [0, 1] \to \mathbb{M}(A)$ as in Proposition 7.3. We choose which $X$-points in each harvestable fitted shell containing 0 or 3 or more points to delete according to $S$, with randomness from the fitted shell's "nearest" 1- and 2-fitted shells, with respect to 1- and 2-Voronoi tessellations. Recall $D : G \times [0, 1] \to \mathbb{M}$ is the independent splitting map on individual process points. We choose which $X$-points outside of fitted shells to delete via $D$, with randomness from each point's nearest 1- and 2-fitted shells. Set $\mathcal{O} := X|_{F(X)^c} \cup F(X) \setminus F'(X)$, so an element of $\mathcal{O}$ is either a fitted shell or $X$-point as described in this paragraph.

Let $V \in \mathcal{V}^*(X)$ with landmark $\ell \in \ell^*(X)$. The 1-Voronoi cell $V$ inherits a $G$-invariant ordering from $Q(\ell^{-1}V) \subseteq [0, 1]$ almost everywhere. Label elements in $X|_{F(X)^c} \cap V$ as $o_{V,k}$ for $1 \leq k \leq N(F(X)^c \cap V)$ so that $Q(\ell^{-1}o_{V,k}) < Q(\ell^{-1}o_{V,k+1})$, and label elements in $F(X) \setminus F'(X) \cap V$ as $o_{V,k}$ for $N(F(X)^c \cap V) < k \leq |\mathcal{O} \cap V|$ in the same way. Similarly, label elements in $\mathcal{O} \cap V^{**}$ as well.

Each 1- or 2-Voronoi cell contains a unique 1- or 2-fitted shell. For $F_i \in F^*(X)$, denote the 1-Voronoi cell containing $F_i$ as $\mathcal{V}(F_i)$, and for $F_j^{**} \in F^{**}(X)$, denote the 2-Voronoi cell containing $F_j^{**}$ as $\mathcal{V}(F_j^{**})$. Let $o \in \mathcal{O}$. We have

$$o = o_{\mathcal{V}(F_i), k} = o_{\mathcal{V}(F_j^{**}), n} =: o(F_i, k; F_j^{**}, n)$$

for some $\mathcal{V}(F_i) \in \mathcal{V}^*(X)$, $\mathcal{V}(F_j^{**}) \in \mathcal{V}^{**}(X)$, $1 \leq k \leq |\mathcal{O} \cap \mathcal{V}(F_i)|$ and $1 \leq n \leq |\mathcal{O} \cap \mathcal{V}(F_j^{**})|$. If $o \in X|_{F(X)^c}$, then $1 \leq k \leq N(F(X)^c \cap \mathcal{V}(F_i))$.

Now we define $\Phi : \mathbb{X}_\alpha \to \mathbb{X}_\beta$ by

$$\Phi(X) := \sum_{F_i^* \in F^*(X)} X|_{F_i^*}$$

$$+ \sum_{F_i \in F^*(X)} \sum_{k=1}^{N(F(X)^c \cap \mathcal{V}(F_i))} D\big(o(F_i, k; F_j^{**}, n), b(\xi(F_i))_k \oplus b(\xi(F_j^{**}))_n\big)$$

$$+ \sum_{F_i \in F^*(X)} \sum_{k=N(F(X)^c \cap \mathcal{V}(F_i))+1}^{|\mathcal{O} \cap \mathcal{V}(F_i)|} S\big(o(F_i, k; F_j^{**}, n), b(\xi(F_i))_k \oplus b(\xi(F_j^{**}))_n\big).$$

We claim $\Phi$ is a monotone, finitary factor from $\mathbb{X}_\alpha$ to $\mathbb{X}_\beta$, and $X - \Phi(X)$ is Poisson point process on $G$ with intensity $\alpha - \beta$.

Proposition 4.1, Remark 7.2 and Proposition 7.3 imply $\Phi(X)$ is a Poisson point process with intensity $\beta$ and $X - \Phi(X)$ is a Poisson process with intensity $\alpha - \beta$. The finitary property follows from Remark 3.3 and Corollary 3.8. Remark 7.2 and Proposition 7.3 imply $\Phi$ is monotone. Since $\xi$ and the $Q$-ordering on Voronoi cells are $G$-invariant, and 1- and 2-Voronoi tessellations are $G$-equivariant, $\Phi$ is $G$-equivariant. □

With coinduction, we extend the results of Theorem 7.4 to prove Theorem 1.4.

PROOF OF THEOREM 1.4. Let $G$ be a nondiscrete, noncompact, locally compact group and let $\alpha > \beta > 0$. If $G$ is compactly generated, Theorem 7.4 implies the result. So, suppose $G$ is not compactly generated, but contains a noncompact, compactly generated open subgroup $H$.

Let $X$ be a Poisson point process on $G$ with intensity $\alpha$, $\psi$ be the map from Lemma 5.4, and $\Phi$ be the map from Theorem 7.4. By Lemma 5.4, $\psi^{-1}(X) = \{X_i\}_{g_i H \in G/H}$ where $X_i$ is a Poisson point process on $H$ with intensity $\alpha$ and $X|_{g_i H} = X_i$ for all $g_i H \in G/H$. Then $\{\Phi(X_i)\}_{g_i H \in G/H}$ is a sequence of Poisson point processes on $H$ with intensity $\beta$, and $\{X_i - \Phi(X_i)\}_{g_i H \in G/H}$ is a sequence of Poisson processes on $H$ with intensity $\alpha - \beta$. Applying $\psi$ to $\{\Phi(X_i)\}_{g_i H \in G/H}$ yields a Poisson point process on $G$ with intensity $\beta$ with the desired properties.

Therefore, the map $\psi \circ \Phi^i \circ \psi^{-1} : \mathbb{X}_\alpha \to \mathbb{X}_\beta$ is a monotone, finitary factor such that $X - \psi(\Phi^i(\psi^{-1}(X)))$ is a Poisson point process on $G$ with intensity $\alpha - \beta$. □

Our construction of the factor and isomorphism maps in the case when $G$ is nondiscrete, noncompact, locally compact, without a noncompact, compactly generated open subgroup relies on Lemma 6.1, and our splitting map in Theorem 1.4 relies on $\lambda(A(e, 6, \rho)) > C_r$ (see Lemma 3.1). We do not know if it is possible to add properties of $S$ in Lemma 6.1 in order to find a constant $\rho$ such that the coupling from Lemma 3.1 could apply, so we do not extend the splitting map to this case.

QUESTION 2. Let $G$ be a nondiscrete, noncompact, locally compact group, without a noncompact, compactly generated open subgroup. Does there exist a splitting factor map for Poisson systems over $G$ with the same properties as the map in Theorem 1.4?

We also note the splitting map in [8] for the $\mathbb{R}^n$ case is equivariant with respect to isometries of $\mathbb{R}^n$, like the factor map in [8] and the isomorphism in [22], in contrast to our $G$-equivariant splitting map. But as mentioned in Section 2, it is known that equivariance cannot extend beyond Isom $\mathbb{R}^n$ for the $\mathbb{R}^n$ case [6].

QUESTION 3. Let $G$ be a nondiscrete, noncompact, locally compact group that is not $\mathbb{R}^n$. If a splitting factor map for Poisson systems over $G$ with the same properties as the map in Theorem 1.4 exists, when does its equivariance extend beyond $G$ itself, and how far can it extend?

**Acknowledgments.** The author thanks Lewis Bowen for helpful comments and discussions provided throughout the writing and revising of this paper, particularly those on word metrics, coinduction and generating sets of compact groups, and Brandon Seward for inspiring the project by suggesting the construction in [22] might work for Lie groups.
The author also thanks the referee for their many facilitative and insightful comments.

**Funding.** This research was supported in part by NSF Grant DMS-1937215.